\documentclass[12pt]{article}
\usepackage{mathrsfs}
\usepackage{epic,eepic,epsf,epsfig}
\usepackage{amsfonts,srcltx,mathrsfs}
\usepackage{multirow}
\usepackage{amsmath}
\usepackage{amssymb}
\usepackage{amsbsy}
\usepackage{graphicx}
\usepackage{amsfonts}%,srcltx}
\usepackage{color}
\usepackage{setspace}
\newtheorem{lem}{Lemma}[section]%
\newtheorem{theorem}[lem]{Theorem}%
\newtheorem{cor}[lem]{Corollary}%

\def\nd{\mathrel{\bigm|\kern-.7em/}}

\def\f{\noindent}

\def\P\GammaL{\hbox{\rm P\GammaL}}

\def\mod{\hbox{\rm mod }}

\begin{document}
\title{How to find all extremal graphs using symmetric subgraphs}

\footnotetext{E-mails: zhangwq@pku.edu.cn}

\author{Wenqian Zhang\\
{\small School of Mathematics and Statistics, Shandong University of Technology}\\
{\small Zibo, Shandong 255000, P.R. China}}
\date{}
\maketitle

\begin{abstract}
 Let $\mathcal{F}$ be a finite family of graphs with $\min_{F\in \mathcal{F}}\chi(F)=r+1\geq3$, where $\chi(F)$ is the chromatic number of $F$.  Set $t=\max_{F\in\mathcal{F}}|F|$.  Let ${\rm EX}(n,\mathcal{F})$ be the set of graphs with maximum edges among all the graphs of order $n$ without any $F\in\mathcal{F}$ as a subgraph. Let $T(n,r)$ be the Tur\'{a}n graph of order $n$ with $r$ parts. Assume that some $F_{0}\subseteq\mathcal{F}$ is a subgraph of the graph obtained from $T(rt,r)$ by embedding a path in its one part.  Simonovits \cite{S1} introduced the concept of symmetric subgraphs, and proved that there exist graphs in ${\rm EX}(n,\mathcal{F})$ which have symmetrical property. In this paper, we aim to find a way to characterize all the extremal graphs for such $\mathcal{F}$ using symmetric subgraphs. Some new extremal results are obtained.

\bigskip

\f {\bf Keywords:} extremal graph; symmetric subgraph; the progressive induction; blow up of graphs.

\end{abstract}

\baselineskip 17 pt

\section{Introduction}
All graphs considered in this paper are finite, undirected and without multi-edges or loops. We first introduce some elementary notations in graph theory. For a graph $G$, let $\overline{G}$ denote its complement. The vertex set and edge set of $G$ are denoted by $V(G)$ and $E(G)$, respectively. Let $|G|=|V(G)|$ and $e(G)=|E(G)|$. For a subset $S\subseteq V(G)$, let $G[S]$ be the subgraph of $G$ induced by $S$, and let $G-S=G[V(G)-S]$. For a vertex $u$, denote $G-\left\{u\right\}$ by $G-u$. Let $N_{G}(u)$ be the set of neighbors of $u$ in $G$, and let $d_{G}(u)=|N_{G}(u)|$. Let $\delta(G)$ denote the minimum degree of $G$. For two disjoint subsets $U,S\subseteq V(G)$, let $e_{G}(U,S)$ be the number of edges between $U$ and $S$. For any terminology used but not defined here, one may refer to \cite{B}.
 
For $\ell\geq2$  graphs $G_{1},G_{2},...,G_{\ell}$, let $\cup_{1\leq i\leq \ell}G_{i}$ (or $G_{1}\cup G_{2}\cup\cdots\cup G_{\ell}$) denote the disjoint union of them. Let $\otimes_{1\leq i\leq \ell}G_{i}$ (or $G_{1}\otimes G_{2}\otimes\cdots\otimes G_{\ell}$) denote the graph obtained from $\cup_{1\leq i\leq \ell}G_{i}$ by connecting each vertex in $G_{i}$ to each vertex in $G_{j}$ for any $1\leq i\neq j\leq \ell$. For a positive integer $t$ and a graph $H$, let $tH=\cup_{1\leq i\leq \ell}H$. For a certain integer $n$, let $K_{n}, C_{n}$ and $P_{n}$ be the complete graph, the cycle and the path of order $n$, respectively. For $r\geq2$, let $K_{n_{1},n_{2},...,n_{r}}$ be the complete $r$-partite graph with parts of sizes $n_{1},n_{2},...,n_{r}$. Let $T(n,r)$ be the Tur\'{a}n graph of order $n$ with $r$ parts (i.e., the complete $r$-partite graph of order $n$, in which each part has $\lfloor\frac{n}{r}\rfloor$ or $\lceil\frac{n}{r}\rceil$ vertices). Assume that $T(n,1)=\overline{K_{n}}$, and $T(n,0)$ is an empty graph. For a graph $G$, let $\chi(G)$ denote its {\em chromatic number}.

Let $F$ be a graph. For a graph $G$, we say $F\subseteq G$ if $G$ contains a copy of $F$, and $G$ is $F$-free otherwise. Let $\mathcal{F}$ be a family of graphs. We say $G$ is $\mathcal{F}$-free, if $G$ is $F$-free for any $F\in\mathcal{F}$. We also uso $F$ to denote $\mathcal{F}$, if $\mathcal{F}=\left\{F\right\}$. Let ${\rm EX}(n,\mathcal{F})$ be the set of $\mathcal{F}$-extremal graphs of order $n$ (i.e., the $\mathcal{F}$-free graphs of order $n$ with maximum edges). Let ${\rm ex}(n,\mathcal{F})$ be the number of edges of any graph in ${\rm EX}(n,\mathcal{F})$.
The classical Tur\'{a}n Theorem \cite{T} states that ${\rm EX}(n,K_{r+1})=\left\{T(n,r)\right\}$,  which is considered as the origin of extremal graph theory. Erd\H{o}s, Stone and Simonovits \cite{E Simonovits,E Stone} proposed the stability theorem
$${\rm ex}(n,\mathcal{F})=(1-\frac{1}{r})\frac{n^{2}}{2}+o(n^{2}),$$ where $\mathcal{F}$ is a finite family of graphs with $\min_{F\in\mathcal{F}}\chi(F)=r+1\geq3$.
The Tur\'{a}n type problem has attracted many researchers (for example, see \cite{CGPW,CY,DJ,EG,ES,EFGG,FS,G,HQL,HQL1,HLF,L,NKSZ,PX,S,S1,S2,WHLM,XKXZ,XZ,Y,Y1,Y2,Y3,ZKS,ZC,ZY,Z}).

Simonovits \cite{S1} introduced the concept of symmetric subgraphs. Let $G$ be a graph. For $\tau\geq1$ induced subgraphs of $G$: $Q_{1},Q_{2},..,Q_{\tau}$, they are called symmetric subgraphs in $G$, if they are connected and vertex-disjoint, and there are isomorphisms $\psi_{j}: Q_{1}\rightarrow Q_{j}$ for all $2\leq j\leq \tau$ such that for any $u\in V(Q_{1})$ and $v\in V(G)-(\cup_{1\leq i\leq\tau}V(Q_{i}))$, $uv\in E(G)$ if and only if $\psi_{j}(u)v\in E(G)$. $u$ and $\psi_{j}(u)$ are called corresponding vertices. These symmetric subgraphs are called symmetric vertices if $|Q_{i}|=1$. Clearly, $\tau$ vertices $v_{1},v_{2},...,v_{\tau}$ of $G$ are symmetric if and only if $\left\{v_{1},v_{2},...,v_{\tau}\right\}$ is independent and $N_{G}(v_{1})=N_{G}(v_{i})$ for any $1\leq i\leq\tau$.

\medskip

\f{\bf Observation 1.} Let $F$ be a graph of order $t\geq2$. Let $G$ be a $F$-free graph. Assume that $Q_{1},Q_{2},...,Q_{t}$ are symmetric subgraphs of $G$. Let $G'$ be the graph obtained from $G$ by adding a new copy of $Q_{1}$, say $Q$, such that $Q,Q_{1},Q_{2},...,Q_{t}$ are symmetric subgraphs in $G'$. Then $G'$ is also $F$-free. In fact, if $F\subseteq G'$, then $V(F)\cap V(Q)\neq\emptyset$ since $G$ is $F$-free. Noting $t=|F|$, there is at least one of $Q_{1},Q_{2},...,Q_{t}$, say $Q_{t}$, such that $V(F)\cap V(Q_{t})=\emptyset$. Let $F'$ be the subgraph of $G$ induced by $(V(F)-V(Q))\cup V(Q_{t})$. By symmetry of $Q$ and $Q_{t}$ in $G'$, we see that $F\subseteq F'\subseteq G$, a contradiction. Hence $G'$ is $F$-free.

\medskip

The above observation is very important in the study of extremal graphs using symmetric subgraphs. For certain integers $n,r$ and $c$, let $\mathbb{D}(n,r,c)$ denote the family of graphs $G$ of order $n$ satisfying the following symmetry condition:\\
$(i)$  It is possible to omit at most $c$ vertices of $G$ so that the remaining graph $G'$ is of form $\otimes_{1\leq i\leq r}G_{i}$, where $||V(G_{i})-\frac{n}{r}|\leq c$ for any $1\leq i\leq r$.\\
$(ii)$ For any $1\leq i\leq r$, $G_{i}=\cup_{1\leq j\leq k_{i}}H^{i}_{j}$, where $H^{i}_{1},H^{i}_{2},...,H^{i}_{k_{i}}$ are symmetric subgraphs in $G$ and $|H^{i}_{1}|\leq c$.

Simonovits \cite{S1} gave the following influential result on extremal graphs. Given integers $r,t$, when $n$ is sufficiently large, $M:=M(r,t)$ is called a constant if it is just a function of $r$ and $t$ (i.e., independent to $n$).

\begin{theorem}{\rm (\cite{S1})}\label{Simonovits 1}
Let $\mathcal{F}$ be a finite family of graphs with $\min_{F\in \mathcal{F}}\chi(F)=r+1\geq3$. Set $t=\max_{F\in\mathcal{F}}|F|$. Assume that $F\subseteq P_{t}\otimes T(t(r-1),r-1)$ for some $F\in\mathcal{F}$. Then, for sufficiently large $n$, $\mathbb{D}(n,r,c)$ contains a graph $G$ in ${\rm EX}(n,\mathcal{F})$, where $c$ is constant. Furthermore, if $G$ is the only extremal graph in $\mathbb{D}(n,r,c)$, then it is
unique in ${\rm EX}(n,\mathcal{F})$.
\end{theorem}

When $\mathbb{D}(n,r,c)$ contains only one $\mathcal{F}$-free graph with maximum edges, then there is only one extremal graph for $\mathcal{F}$ by Theorem \ref{Simonovits 1}. 
 We should also note that $\mathbb{D}(n,r,c)$ usually contains many extremal graphs. If so, Theorem \ref{Simonovits 1} will not help us to find all extremal graphs. Motivated by this, we seek a way to characterize all extremal graphs in this paper.  The following lemma is one main result of this paper.  We will use it together with the progressive induction \cite{S} to obtain some new extremal results.

\begin{lem}\label{Pl-free}
Let $\mathcal{F}$ be a finite family of graphs with $\min_{F\in \mathcal{F}}\chi(F)=r+1\geq3$. Set $t=\max_{F\in\mathcal{F}}|F|$. Assume that $F_{1}\subseteq P_{t}\otimes T(t(r-1),r-1)$ and $F_{2}\subseteq tK_{1,t}\otimes tK_{1,t}\otimes T(t(r-2),r-2)$ for some $F_{1},F_{2}\in\mathcal{F}$. For any given constant $N_{0}:=N_{0}(r,t)$, there exists a constant $O(r,t)$ such that any graph $G$ in ${\rm EX}(n,\mathcal{F})$ with $n$ sufficiently large, has an induced subgraph $\otimes_{1\leq i\leq r}G_{i}$ satisfying:\\
 $(i)$  $|G_{i}|=O(r,t)\cdot N_{0}$ for any $1\leq i\leq r$;\\
 $(ii)$ for any $1\leq i\leq r$, $G_{i}=\cup _{1\leq j\leq \ell_{i}}H^{i}_{j}$, where $H^{i}_{1},H^{i}_{2},...,H^{i}_{\ell_{i}}$ are symmetric subgraphs of $G$ and $|H^{i}_{1}|\leq O(r,t)$.
\end{lem}

For a finite family of graphs $\mathcal{F}$ with $\min_{F\in \mathcal{F}}\chi(F)=r+1\geq3$ and $t=\max_{F\in\mathcal{F}}|F|$, define $q(\mathcal{F})$ to be the minimum integer $s$ such that $F\subseteq \overline{K_{s}}\otimes T(tr,r)$ for some $F\subseteq\mathcal{F}$. Clearly, $1\leq q(\mathcal{F})\leq t$. 

\begin{theorem}\label{matching-free}
Let $\mathcal{F}$ be a finite family of graphs with $\min_{F\in \mathcal{F}}\chi(F)=r+1\geq3$. Set $t=\max_{F\in\mathcal{F}}|F|$ and $q=q(\mathcal{F})$. Assume that $F_{1}\subseteq tK_{2}\otimes T(t(r-1),r-1)$ and $F_{2}\subseteq (K_{q-1,t}\cup K_{1,t})\otimes T(t(r-1),r-1)$ for some $F_{1},F_{2}\in\mathcal{F}$. Then for sufficiently large $n$, any graph $G$ in ${\rm EX}(n,\mathcal{F})$ has a partition $V(G)=W\cup(\cup_{1\leq i\leq r}S_{i})$ satisfying:\\
$(i)$ $|W|=q-1$ and $|S_{i}|=\lfloor\frac{n-q+1}{r}\rfloor$ or $\lceil\frac{n-q+1}{r}\rceil$ for any $1\leq i\leq r$;\\
$(ii)$ for any $1\leq i\leq r$, there exists $S'_{i}\subseteq S_{i}$ such that $|S_{i}-S'_{i}|\leq t^{2}$, and $N_{G}(v)=V(G)-S_{i}$ for any $v\in S'_{i}$.
\end{theorem}

Clearly, all the extremal graphs in Theorem \ref{matching-free} are contained in $\mathbb{D}(n,r,c)$ with $c=q+rt^{2}$. Note that in Theorem \ref{matching-free}, the condition  $F_{2}\subseteq (K_{q-1,t}\cup K_{1,t})\otimes T(t(r-1),r-1)$ for some $F_{2}\in\mathcal{F}$  is not removable. For example, let $\mathcal{F}=\left\{3K_{2}\otimes\overline{K_{10}},P_{5}\otimes\overline{K_{10}}\right\}$.
When $n$ is large, we can show that $G=H\otimes\overline{K_{\lfloor\frac{n}{2}\rfloor}}$ is an extremal graph for $\mathcal{F}$, where $H$ is any double star. Clearly, when the two center vertices have almost equal degree in $H$, the resulting extremal graph $G$ is not in  $\mathbb{D}(n,r,c)$ for any constant $c$.

Simonovits (communicated with Liu \cite{L}) asked the following question:
 
 \medskip
 
 \f{\bf Question 1.} Characterize graphs in $\mathcal{F}$ with $\min_{F\in \mathcal{F}}\chi(F)=r+1\geq3$, such that $K_{q-1}\otimes T(n-q+1,r)$  is the unique extremal graph.

\medskip

Let $I^{12}$ denote the icosahedron. Simonovits \cite{S2} observed that $I^{12}$ has no independent set of order $4$. Thus, $\chi(I^{12})=4$ and $K_{2}\otimes T(n-2,3)$ is $I^{12}$-free. Moreover, $I^{12}\subseteq F_{i}$ for any $1\leq i\leq3$, where $F_{i}$ is defined as in Theorem \ref{Icosahe}. As a classical result, Simonovits \cite{S2} characterized the unique extremal graph of large order for $I^{12}$, by showing the following stronger result.

\begin{theorem}{\rm (\cite{S2})} \label{Icosahe}
Let $\mathcal{F}=\left\{F_{1},F_{2},F_{3}\right\}$, where
$$F_{1}=P_{6}\otimes 3K_{1}\otimes3K_{1},$$
$$F_{2}=(K_{1,2}\cup K_{2})\otimes(K_{2}\cup 2K_{1})\otimes3K_{1},$$
$$F_{3}=2K_{3}\otimes3K_{1}\otimes3K_{1}.$$
Then ${\rm EX}(n,\mathcal{F})=\left\{K_{2}\otimes T(n-2,3)\right\}$ for sufficiently large $n$. 
Consequently, ${\rm EX}(n,I^{12})=\left\{K_{2}\otimes T(n-2,3)\right\}$ for sufficiently large $n$.
 \end{theorem}
 
\medskip

At the end of \cite{S2}, Simonovits also hopes a generalized result. In this paper, we generalize Theorem \ref{Icosahe} (by letting $\ell=6,r=3$) as follows, which supplies more examples for Question 1. 

For certain integers $n,r$, let $K^{+}_{n_{1},n_{2},...,n_{r}}$ be the graph of order $n$ obtained from the complete $r$-partite graph $K_{n_{1},n_{2},...,n_{r}}$ by adding a (near) perfect matching in each of its $r$ parts, where $\sum_{1\leq i\leq r}n_{i}=n$ and $|n_{i}-n_{j}|\leq2$ for any $1\leq i<j\leq r$. Let $\mathcal{G}_{n,r}$ be the set of all such graphs $K^{+}_{n_{1},n_{2},...,n_{r}}$. Let $G_{n,r}$ be the graph obtained from $T(n,r)$ by adding a (near) perfect matching in each of its $r$ parts. Clearly, $G_{n,r}$ is an (but not only) extremal graph with maximum edges among the graphs in  $\mathcal{G}_{n,r}$.

\begin{theorem}\label{general Icosahe}
Assume that $\ell\geq4$ is an even integer, and $r,m$ are two integers with $2\leq r\leq\ell-2$ and $m\geq\ell+1$. Let $\mathcal{F}=\left\{F_{1},F_{2},F_{3}\right\}$, where
$$F_{1}=(P_{\ell}\cup\overline{K_{m}})\otimes T(m(r-1),r-1),$$
$$F_{2}=(K_{1,2}\cup mK_{2})\otimes mK_{2}\otimes T(m(r-2),r-2),$$ $$F_{3}=(K_{\frac{\ell}{2}}\cup K_{\ell-1}\cup\overline{K_{m}})\otimes T(m(r-1),r-1).$$
Then  for sufficiently large $n$,\\
$(i)$. ${\rm EX}(n,\mathcal{F})=\left\{K_{\frac{\ell-2}{2}}\otimes T(n-\frac{\ell-2}{2},r)\right\}$ for $2\leq r\leq\ell-3$;\\
$(ii)$. ${\rm EX}(n,\mathcal{F})\subseteq\left\{K_{\frac{r}{2}}\otimes T(n-\frac{r}{2},r)\right\}\cup \mathcal{G}_{n,r}$ for $2\leq r=\ell-2$.
\end{theorem}

\begin{theorem}\label{m Icosahe}
Assume that $m\geq2$ is an integer.
Then  ${\rm EX}(n,m I^{12})=\left\{K_{3m-1}\otimes T(n-3m+1,3)\right\}$ for sufficiently large $n$.
\end{theorem}

The following Theorem \ref{bi-matching-free} generalizes the extremal results for the $p$-th power of paths \cite{XKXZ,Y1} and the odd prism \cite{HLF}.

\begin{theorem}\label{bi-matching-free}
For $\ell\geq3$ and $m\geq2\ell$, let $\mathcal{F}=\left\{F_{1},F_{2},F_{3},F_{4}\right\}$, where
$$F_{1}=(P_{\ell}\cup\overline{K_{m}})\otimes T(m(r-1),r-1),$$
$$F_{2}=mK_{2}\otimes mK_{2}\otimes T(m(r-2),r-2),$$ 
$$F_{3}=(2K_{\ell-1}\cup\overline{K_{m}})\otimes (K_{2}\cup\overline{K_{m}})\otimes T(m(r-2),r-2),$$ 
$$F_{4}=(K_{\ell-1}\cup\overline{K_{m}})\otimes (K_{\ell-1}\cup\overline{K_{m}})\otimes T(m(r-2),r-2).$$ 
 For sufficiently large $n$, any graph $G\in{\rm EX}(n,\mathcal{F})$ satisfies
 $G=H\otimes T(n-|H|,r-1)$, where $H$ is an extremal $P_{\ell}$-free graph.
\end{theorem}

The rest of this paper is organized as follows. In Section 2, we apply our results on the extremal problems for blow up or odd-ballooning of graphs. In Section 3, we include some lemmas needed in the proofs of the main results. In Section 4, we give the proof of Lemma \ref{Pl-free}. In Section 5, we give the proof of Theorem \ref{matching-free}. In Section 6, we give the proofs of Theorem \ref{general Icosahe} and Theorem \ref{m Icosahe}. The proof of Theorem \ref{bi-matching-free} is given in Section 7.

\section{Applications on extremal graphs for blow up or odd-ballooning of graphs}

Let $\mathcal{F}$ be a finite family of graphs with $\min_{F\in\mathcal{F}}\chi(F)=r+1\geq3$. Let
$\mathcal{M}(\mathcal{F})$ be the family of minimal graphs (with no isolated vertices) $M$ satisfying the following: there exist an $F\in\mathcal{F}$ such that $F\subseteq(M\cup\overline{K_{t}})\otimes T((r-1)t,r-1)$, where $t=|F|$. We call $\mathcal{M}(\mathcal{F})$ the decomposition family of $F$. The
decomposition family $\mathcal{M}(\mathcal{F})$ always contains some bipartite graphs.

A covering of a graph is a set of vertices incident with all edges of the graph.
An independent covering of a bipartite graph is an independent set incident with all edges.
Let $\beta(F)$ denote the covering number
of $F$, i.e., the minimum number of vertices in a covering of $F$. If $F$ is bipartite, let $q(F)$ be the minimum order of an independent covering of $F$. 
Let $\mathcal{F}$ be a finite family of graphs containing at least one bipartite graph. 
Define $q(\mathcal{F})=\min\left\{q(F)~|~F\in\mathcal{F}~is~bipartite\right\}$. This definition of $q(\mathcal{F})$ is consistent with the one in the introduction.
The subgraph covering family $\mathcal{B}(\mathcal{F})$ of $\mathcal{F}$ is the set of subgraphs $F[U]$ of $F\in\mathcal{F}$, where $U$ is a covering of $F$ with $|U|\leq q(\mathcal{F})-1$. If $\beta(F)\geq q(F)$ for each $F\in\mathcal{F}$, set $\mathcal{B}(\mathcal{F})=\left\{K_{q(\mathcal{F})}\right\}$.

The following result is an easy consequence of Theorem \ref{matching-free}. (There is no edge inside $S_{i}$ in the condition that $F_{2}\subseteq K_{q-1,t}\cup K_{2}$ with some $F_{2}\in\mathcal{M}(\mathcal{F})$.)

\begin{cor}\label{matching one degree}
Let $\mathcal{F}$ be a finite family of graphs with $\min_{F\in \mathcal{F}}\chi(F)=r+1\geq3$. Set $t=\max_{F\in\mathcal{F}}|F|$ and $q=q(\mathcal{F})$. Assume that $F_{1}\subseteq tK_{2},F_{2}\subseteq K_{q-1,t}\cup K_{2}$ for some $F_{1},F_{2}\in\mathcal{M}(\mathcal{F})$. Then for sufficiently large $n$, any graph $G\in{\rm EX}(n,\mathcal{F})$ is of form
 $G=H\otimes T(n-q+1,r)$, where $H\in{\rm EX}(q-1,\mathcal{B}(\mathcal{F}))$.
\end{cor}

Corollary \ref{matching one degree} has concrete examples (see \cite{CY,Y2}). Given a graph $G$ and an integer $p$, the blow-up of $G$, denoted by $G^{p+1}$, is the
graph obtained from  $G$ by replacing each edge in $G$ by a clique of size $p+1$, where the new
vertices of the cliques are all distinct. Liu \cite{L} proved the following result.

\begin{theorem}{\rm (\cite{L})}
For any $p\geq2$ and $k\geq3$, when $n$ is sufficiently large, we have the following
results:\\
$(i)$ For any $C^{p+1}_{k}\neq C^{3}_{3}$, the unique graph in ${\rm EX}(n,C^{p+1}_{k})$ is $K_{\lfloor\frac{k-1}{2}\rfloor}\otimes T(n-\lfloor\frac{k-1}{2}\rfloor,p)$ for odd $k$; the unique graph in ${\rm EX}(n,C^{p+1}_{k})$ is obtained from $K_{\lfloor\frac{k-1}{2}\rfloor}\otimes T(n-\lfloor\frac{k-1}{2}\rfloor,p)$ by adding one edge in one part of $T(n-\lfloor\frac{k-1}{2}\rfloor,p)$ for even $k$.\\
$(ii)$ For $C^{3}_{3}$, if $4|n$, the graph $G^{*}_{n}$ (obtained from $T(n,2)$ by adding a (near) perfect matching in each of its two parts) is the unique graph in ${\rm EX}(n,C^{3}_{3})$; otherwise both $G^{*}_{n}$ and $K_{1}\otimes T(n-1,2)$ are in ${\rm EX}(n,C^{3}_{3})$.
\end{theorem}

For the extremal graphs for $C^{3}_{3}$, one can see that 
$$C^{3}_{3}\subseteq P_{4}\otimes\overline{K_{2}},K_{1,2}\otimes(K_{2}\cup K_{1}),K_{3}\otimes\overline{K_{3}}.$$ 
From $(ii)$ of Theorem \ref{general Icosahe} (by letting $r=2,\ell=4$), we can obtain ${\rm EX}(n,C^{3}_{3})\subseteq\left\{K_{1}\otimes T(n-1,2)\right\}\cup\mathcal{G}_{n,2}$. Thus, we have the following result (after a simple calculation).

\begin{cor}
For $C^{3}_{3}$ and large $n$, we have\\
${\rm EX}(n,C^{3}_{3})=\left\{K^{+}_{\frac{n}{2},\frac{n}{2}}\right\}$ for $n\equiv0~(\mod4)$; \\
${\rm EX}(n,C^{3}_{3})=\left\{K^{+}_{\frac{n-1}{2},\frac{n+1}{2}},K_{1}\otimes T(n-1,2)\right\}$ for $n\equiv1,3~(\mod4)$;\\
 ${\rm EX}(n,C^{3}_{3})=\left\{K^{+}_{\frac{n}{2},\frac{n}{2}},K^{+}_{\frac{n-2}{2},\frac{n+2}{2}},
 K_{1}\otimes T(n-1,2)\right\}$ for $n\equiv2~(\mod4)$.
\end{cor}

For a general graph $F$, the extremal graphs for $F^{p+1}$ with $p\geq\chi(F)+1$ have a specified structure as follows (using Theorem \ref{matching-free} and a characterization for $\mathcal{M}(F^{p+1})$ in \cite{L}).

\begin{theorem}\label{blow up}
Let $F$ be a graph, and let $p\geq\chi(F)+1$ be an integer. Set $q=q(\mathcal{M}(F^{p+1}))$. Then for sufficiently large $n$, any graph $G$ in ${\rm EX}(n,F^{p+1})$ has a partition $V(G)=W\cup(\cup_{1\leq i\leq p}S_{i})$ satisfying:\\
$(i)$ $|W|=q-1$ and $|S_{i}|=\lfloor\frac{n-q+1}{p}\rfloor$ or $\lceil\frac{n-q+1}{p}\rceil$ for any $1\leq i\leq p$;\\
$(ii)$ there is a constant $M:=M(|F|,p)$, such that for any $1\leq i\leq p$, there exists $S'_{i}\subseteq S_{i}$ such that $|S_{i}-S'_{i}|\leq M$, and $N_{G}(v)=V(G)-S_{i}$ for any $v\in S'_{i}$.
\end{theorem}

Note that in Theorem \ref{blow up}, when $F$ is non-bipartite, the graphs in ${\rm EX}(n,F^{p+1})$ (with $p\geq\chi(F)+1$) have been completely characterized by Yuan \cite{Y2}.  The extremal graphs for $F^{p+1}$ are also characterized when $F$ is a tree (see \cite{CY,G,L,WHLM}) and other specified graphs (\cite{NKSZ}).

For a graph $F$, the odd-ballooning of $F$ is a graph obtained from $F$, such that each edge in $F$ is replaced by an odd cycle and all new vertices for different odd cycles are distinct. It seems that $F$ is more concerned when it is a bipartite graph (see \cite{HQL,HQL1,Y3,ZKS,ZC,ZY}). Let $F=F[A,B]$ with $|A|\leq|B|$ be a bipartite graph. An odd-ballooning of $F$ is called good, if the replacing triangles can exist only when the corresponding edges are incident with a leaf in $B$. Let $F^{good}$ be a good odd-ballooning of $F$. One can check that $\mathcal{M}(F^{good})$ satisfies the condition in Theorem \ref{matching-free}. Thus, we also have a specified characterization.

\begin{theorem}\label{odd-ballooning}
Let  $F=F[A,B]$ with $|A|\leq|B|$ be a bipartite graph, and let $F^{good}$ be a good odd-ballooning of it. Set $q=q(\mathcal{M}(F^{odd}))$. Then for sufficiently large $n$, any graph $G$ in ${\rm EX}(n,F^{odd})$ has a partition $V(G)=W\cup S_{1}\cup S_{2}$ satisfying:\\
$(i)$ $|W|=q-1$ and $|S_{i}|=\lfloor\frac{n-q+1}{2}\rfloor$ or $\lceil\frac{n-q+1}{2}\rceil$ for any $1\leq i\leq 2$;\\
$(ii)$ there is a constant $M:=M(|F^{odd}|)$, such that for any $1\leq i\leq 2$, there exists $S'_{i}\subseteq S_{i}$ such that $|S_{i}-S'_{i}|\leq M$, and $N_{G}(v)=V(G)-S_{i}$ for any $v\in S'_{i}$.
\end{theorem}

Let $T=T[A,B]$ with $|A|\leq|B|$ be a tree. Define $\delta(A)=\min_{u\in A}d_{T}(u)$. Let $T_{0}$ be a good odd-ballooning of $T$. Let $\delta_{1}(A)$ be the minimum number of triangles containing a vertex $v\in A$ satisfying $d_{T}(v)=\delta(A)$. Using Theorem \ref{Simonovits 1}, Zhu and Chen \cite{ZC} obtained ${\rm ex}(n,T_{0})$ and an extremal graph for any tree $T$. When $\delta_{1}(A)=\delta(A)$, then $T_{0}$ is just the friendship graph (see \cite{EFGG}). When $\delta_{1}(A)<\delta(A)=1$, all the extremal graphs can be obtained by Corollary \ref{matching one degree}. When $\delta_{1}(A)<\delta(A)$ and $\delta(A)\geq2$, their result \cite{ZC} is the following. For certain integers $n,a,k$, let $G_{n,a,k}$ be the graph obtained from $K_{a-1}\otimes T(n-a+1,2)$ by adding a $K_{k-1,k-1}$ in one part of $T(n-a+1,2)$.

\begin{theorem}{\rm (\cite{ZC})}
Let $T=T[A,B]$ be a tree with $a=|A|\leq |B|$ and $\delta(A)=k\geq2$, and let $T_{0}$ be a good odd-ballooning of $T$ such that $\delta_{1}(A)<\delta(A)$.  If $n$ is sufficiently large, then
$${\rm ex}(n,T_{0})=\binom{a-1}{2}+(a-1)(n-a+1)+e(T(n-a+1,2))+(k-1)^{2}.$$
Furthermore, $G_{n,a,k}$ is an extremal graph.
\end{theorem}

Let $G'_{n,a,4}$ be the graph obtained from $K_{a-1}\otimes T(n-a+1,2)$ by adding a $3K_{3}$ in one part of $T(n-a+1,2)$. Using Theorem \ref{odd-ballooning} and a similar discussion as \cite{ZC}, we can obtain the following result.

\begin{cor}
Let $T=T[A,B]$ be a tree with $a=|A|\leq |B|$ and $\delta(A)=k\geq2$, and let $T_{0}$ be a good odd-ballooning of $T$ such that $\delta_{1}(A)=k_{1}<k$.  If $n$ is sufficiently large, then
${\rm EX}(n,T_{0})=\left\{G_{n,a,k}\right\}$ for $(k,k_{1})\neq(4,3)$, and
${\rm EX}(n,T_{0})=\left\{G_{n,a,4},G'_{n,a,4}\right\}$ for $(k,k_{1})=(4,3)$.
\end{cor}

For $2\leq s\leq t$, let $H$ be the graph obtained from $K_{t-1,t-1}$ by deleting a matching of $t-2$ edges. Let $G_{s,t}$ be the graph obtained from $K_{s-1}\otimes T(n-s+1,2)$ by
 embedding $H$ into one class of $T(n-s+1,2)$. Let $G'_{3,3}$
be the graph obtained from $K_{2}\otimes T(n-2,2)$ by
embedding a triangle into one class of $T(n-2,2)$. Using Theorem \ref{Simonovits 1}, Peng and Xia \cite{PX} obtained the following result.

\begin{theorem}{\rm (\cite{PX})}
Let $F_{s,t}$ be the odd‐ballooning of $K_{s,t}$ with $t\geq s\geq2$, where $s+t\geq6$ and each odd cycle has length at least 5. Then for sufficiently large $n$,
$${\rm ex}(n,F_{s,t})=\lceil\frac{n-s+1}{2}\rceil\lfloor\frac{n-s+1}{2}\rfloor+(s-1)(n-s+1)
+\binom{s-1}{2}+t^{2}-3t+3.$$
 Moreover, ${\rm EX}(n,F_{s,t})=\left\{G_{s,t}\right\}$ for $t\geq4$, and ${\rm EX}(n,F_{3,3})\supseteq\left\{G_{3,3},G'_{3,3}\right\}$ for $t=3$.
\end{theorem}

Let $G''_{3,3}$ be the graph obtained from $\overline{K_{2}}\otimes T(n-2,2)$ by
embedding a $C_{4}$ into one class of $T(n-2,2)$. Note that $F_{3,3}$ contains a $F_{1,3}$ after deleting any two non-adjacent vertices. We see that $G''_{3,3}$ is $F_{3,3}$-free (using \cite{HQL1} or \cite{Y3}). (If we also consider the case that $w_{1}$ and $w_{2}$ is not adjacent in Claim 12 of \cite{PX}, we can find $G''_{3,3}$.) Using Theorem \ref{odd-ballooning} and a similar discussion as \cite{PX}, we can obtain the following result.

\begin{cor}
Let $F_{3,3}$ be the odd-ballooning of $K_{3,3}$, such that each odd cycle has length at least 5. Then for sufficiently large $n$, ${\rm EX}(n,F_{3,3})=\left\{G_{3,3},G'_{3,3},G''_{3,3}\right\}$.
\end{cor}

\section{Preliminaries}

To prove the main results of this paper, we first include several lemmas.
The following result is the classical stability theorem proved by  Erd\H{o}s,
Stone and Simonovits \cite{E Simonovits,E Stone}:

\begin{theorem} {\rm (\cite{E Simonovits,E Stone})}\label{stability}
Let $\mathcal{F}$ be a finite family of graphs with $\min_{F\in\mathcal{F}}\chi(F)=r+1\geq3$. For every $\epsilon>0$, there exist $\delta>0$ and $n_{0}$ such that, if $G$ is an $\mathcal{F}$-free graph of order $n\geq n_{0}$ with $e(G)>(\frac{r-1}{2r}-\delta)n^{2}$, then $G$ can be
obtained from $T(n,r)$ by adding and deleting at most $\epsilon n^{2}$ edges.
\end{theorem}

The following well-known lemma is proved by Erd\H{o}s \cite{E} and
Simonovits \cite{S}.

\begin{lem}{\rm (\cite{E,S})}\label{mini degree}
Let $\mathcal{F}$ be a finite family of graphs with $\min_{F\in\mathcal{F}}\chi(F)=r+1\geq3$. For any given $\theta>0$, if $n$ is sufficiently large, then $\delta(G)>(1-\frac{1}{r}-\theta)n$ for any $G\in{\rm EX}(n,\mathcal{F})$.
\end{lem}

\medskip

As one of the earliest extremal
results in graph theory, Erd\H{o}s and Gallai \cite{EG} proved the following theorem. 

\begin{lem}{\rm (\cite{EG})}\label{ex n Pt}
For any $\ell\geq2$ and $n\geq \ell$, let $G$ be a $P_{\ell}$-free graph of order $n$. Then $e(G)\leq\frac{\ell-2}{2}n$ with equality if and only if $n=(\ell-1)t$ and $G=\cup_{1\leq i\leq t}K_{\ell-1}$. 
\end{lem}

Faudree and Schelp \cite{FS} improved this result as follows:

\begin{lem}{\rm (\cite{FS})}\label{ex n Pt sharp}
Assume that $\ell\geq2$ and $n=(\ell-1)t+s$, where $0\leq s<\ell-1$. let $G$ be a $P_{\ell}$-free graph of order $n$. Then $e(G)\leq t\binom{\ell-1}{2}+\binom{s}{2}$ with equality if and only if either $G=(\cup_{1\leq i\leq t}K_{\ell-1})\cup K_{s}$, or $G=(\cup_{1\leq i\leq t-t_{0}}K_{\ell-1})\cup(K_{\frac{\ell}{2}-1}\vee\overline{K_{(\ell-1)t_{0}-\frac{\ell}{2}+s+1}})$ for some $1\leq t_{0}\leq t$ when $\ell$ is even and $s=\frac{\ell}{2}\pm1$. 
\end{lem}

 Simonovits \cite{S} introduced the progressive induction.
 
\medskip

\begin{lem}{\rm (\cite{S})}\label{progressive induction}
Let $\mathfrak{U}=\cup_{1\leq n\leq \infty}\mathfrak{U}_{n}$ be a set of given elements, where $\mathfrak{U}_{n}$ are disjoint finite subsets of $\mathfrak{U}$. Let $B$ be a condition or property defined on $\mathfrak{U}$ (i.e., the elements of $\mathfrak{U}$ may satisfy or
not satisfy $B$). Let $\phi(a)$  be a non-negative integer-valued function defined on $\mathfrak{U}$, such that \\
$(i)$ if $a$ satisfies $B$, then $\phi(a)=0$;\\
$(ii)$ There is a $M_{0}$, such that for any $n>M_{0}$ and $a\in\mathfrak{U}_{n}$, either $a$ satisfies $B$, or there exist an $n'$ with $\frac{n}{2}<n'<n$ and an $a'\in\mathfrak{U}_{n'}$ such that $\phi(a)<\phi(a')$.\\
Then there exists an $N_{0}$ such that any $a\in \mathfrak{U}_{n}$ with $n>N_{0}$ satisfies $B$.
\end{lem}

The following result is a consequence of Lemma 3.2.1 in \cite{S1}. Furthermore, Simonovits pointed out that the orders of the symmetric subgraphs are bounded (i.e., $O(1)$) at the end of the proof of Lemma 3.2.1 in \cite{S1}.

\begin{lem}{\rm (\cite{S1})}\label{symmetric}
Let $G$ be a $P_{\ell}$-free graph of order $n$. For any sequence of positive integers $\pi=\left\{m_{k}\right\}_{k\geq1}$, there exist two positive integers $C(\ell,\pi)$ and $D(\ell,\pi)$ such that if $n\geq C(\ell,\tau)$, then $G$ has $m_{k}$ symmetric subgraphs $Q_{1},Q_{2},...,Q_{m_{k}}$, where $1\leq|Q_{1}|=k\leq D(\ell,\pi)$.
\end{lem}

The following result will be used in our proofs.

\begin{lem}\label{desired symmetric subgraphs}
For $\ell\geq3$, let $G$ be a $P_{\ell}$-free graph. Suppose that $S\subseteq V(G)$ satisfies that $\frac{|S|}{\ell-1}$ is an integer, and $G[S]=\cup_{1\leq i\leq\tau}Q_{i}$, where $Q_{1},Q_{2},...,Q_{\tau}$ are symmetric subgraphs of $G$ with $\tau\geq\ell$. If $e(G[S])+e_{G}(S,V(G)-S)\geq\frac{\ell-2}{2}|S|$, then $e(G[S])+e_{G}(S,V(G)-S)=\frac{\ell-2}{2}|S|$, and we have the following conclusions.\\
$(i)$ For odd $\ell\geq3$, $Q_{i}=K_{\ell-1}$ for $1\leq i\leq\tau$ and $e_{G}(S,V(G)-S)=0$.\\
$(ii)$ For even $\ell\geq4$, either $Q_{i}=K_{\ell-1}$ for any $1\leq i\leq\tau$ and $e_{G}(S,V(G)-S)=0$, or $Q_{i}=K_{1}$ for any $1\leq i\leq\tau$ and there are $\frac{\ell-2}{2}$ vertices in $V(G)-S$ which are all adjacent to the (unique) vertex in $Q_{i}$ for all $1\leq i\leq\tau$.
\end{lem}

\f{\bf Proof:} Let $B\subseteq V(G)-S$ be the set of vertices which have at least one neighbor in $Q_{1}$. Let $H$ be the graph obtained from $G[S\cup B]$ by adding some edges inside $B$ such that $H[B]$ is a complete graph. Let $b=|B|$. If $b=0$, then $e(G[S])\geq\frac{\ell-2}{2}|S|$. Note that
$G[S]$ is $P_{\ell}$-free. By Lemma \ref{ex n Pt}, we must have $e(G[S])=\frac{\ell-2}{2}|S|$ and $Q_{i}=K_{\ell-1}$ for any $1\leq i\leq\tau$. Now we assume $b\geq1$. 
Clearly, $|H|=|S|+b\geq\ell$ and $e(H)=e(G[S])+e_{G}(S,V(G)-S)+\binom{b}{2}\geq\frac{\ell-2}{2}|S|+\binom{b}{2}$.

By Symmetry, each vertex in $B$ has at least one neighbor in $Q_{i}$ for any $1\leq i\leq\tau$. Since $b\geq1$, $H$ is connected. 
Now we show that $H$ is $P_{\ell}$-free. Suppose not, i.e., there is a $P_{\ell}$ in $H$, say $P=u_{1}u_{2}\cdots u_{\ell}$. We choose such a $P$ that $P$ contains as small number of edges in $H[B]$ as possible. If $P$ contains no edges inside $B$, then $P$ is also a path in $G$, a contradiction. Hence, we can assume that $P$ contains at least one edge in $H[B]$. Without loss of generality, assume that $u_{2}u_{3}$ is an edge in $H[B]$.  Noting $\tau\geq\ell$, we can choose one symmetric subgraph, say $Q_{\tau}$, such that $P$ contains no vertices in $V(Q_{\tau})$. 
Since $u_{2},u_{3}\in B$ have neighbors in $Q_{\tau}$, let $v,w\in V(Q_{\tau})$ be neighbors of $u_{2},u_{3}$, respectively. Let $P(v,w)$ be a path from $v$ to $w$ in $Q_{\tau}$. Then $u_{1}u_{2}P(v,w)u_{3}\cdots u_{\ell}$ is a path with less number of edges inside $B$ than $P$. Clearly, its sub-path of order $\ell$ is a $P_{\ell}$ with less number of edges inside $B$ than $P$. This contradicts the choice of $P$. Hence, $H$ is $P_{\ell}$-free. Recall that $|H|=|S|+b$ and $e(H)\geq\frac{\ell-2}{2}|S|+\binom{b}{2}$. Since $H$ is connected and $\frac{|S|}{\ell-1}$ is an integer, by Lemma \ref{ex n Pt sharp}, we must have that $\ell$ is even and $H=K_{\frac{\ell-2}{2}}\otimes\overline{K_{|S|+b-\frac{\ell-2}{2}}}$. Note that $H$ will contain no edges after deleting $\frac{\ell-2}{2}$ vertices in the part $V(K_{\frac{\ell-2}{2}})$.
Since $\tau\geq\ell$, we must have $|Q_{i}|=1$ for any $1\leq i\leq\tau$.  Moreover, $|B|=\frac{\ell-2}{2}$ and each vertex in $B$ is adjacent to the unique vertex in  $|Q_{i}|=1$ for any $1\leq i\leq\tau$.
This completes the proof.
 \hfill$\Box$

\medskip

\section{Proof of Lemma \ref{Pl-free}}

The following fact will be used.

\medskip

\f{\bf Fact 1.}
Let $A_{1},A_{2},...,A_{m}$ be $m\geq2$ finite sets. Then
$$|\cap_{1\leq i\leq m}A_{i}|\geq(\sum_{1\leq i\leq m}|A_{i}|)-(m-1)|\cup_{1\leq i\leq m}A_{i}|.$$

\medskip

\f{\bf Proof of Lemma \ref{Pl-free}.} Recall that $\min_{F\in \mathcal{F}}\chi(F)=r+1\geq3$ and $t=\max_{F\in\mathcal{F}}|F|$. There are  graphs $F_{1},F_{2}\in\mathcal{F}$ such that
$$F_{1}\subseteq P_{t}\otimes T(t(r-1),r-1)$$
 and
$$F_{2}\subseteq tK_{1,t}\otimes tK_{1,t}\prod T(t(r-2),r-2).$$
For sufficiently large $n$, let $G\in {\rm EX}(n,\mathcal{F})$. Let $\theta:=\theta(r,t)>0$ be a small constant (such that all the following inequalities on it are satisfied). By Lemma \ref{mini degree}, we have that $\delta(G)>(\frac{r-1}{r}-\theta)n$.  
As is well known, $e(T(n,r))\geq\frac{r-1}{2r}n^{2}-\frac{r}{8}$. Since $T(n,r)$ is $\mathcal{F}$-free, we have $e(G)\geq e(T(n,r))\geq\frac{r-1}{2r}n^{2}-\frac{r}{8}$.
Since $G$ is $\mathcal{F}$-free and $n$ is large, by Theorem \ref{stability} (letting $\epsilon=\theta^{3}$), we see that $G$ can be obtained from $T(n,r)$ by deleting and adding at most $\theta^{3}n^{2}$ edges. It follows that
 $$(\frac{r-1}{2r}-\theta^{3})n^{2}\leq e(G)\leq(\frac{r-1}{2r}+\theta^{3})n^{2}.$$
  Let $V(G)=\cup_{1\leq i\leq r}V_{i}$ be a partition such that $\sum_{1\leq i\leq r}e(G[V_{i}])$ is minimum. Recall that $G$ can be obtained from $T(n,r)$ by deleting at most $\theta^{3}n^{2}$ edges. Thus, $$\sum_{1\leq i\leq r}e(G[V_{i}])\leq\theta^{3}n^{2}.$$

\medskip

\f{\bf Claim 1.} $||V_{i}|-\frac{n}{r}|\leq\theta n$ for any $1\leq i\leq r$.

\medskip

\f{\bf Proof of Claim 1.} Let $a=\max_{1\leq i\leq r}||V_{i}|-\frac{n}{r}|$. Without loss of generality, assume that $a=||V_{1}|-\frac{n}{r}|$. Using the Cauchy-Schwarz inequality, we see
$$2\sum_{2\leq i<j\leq r}|V_{i}||V_{j}|=(\sum_{2\leq i\leq r}|V_{i}|)^{2}-\sum_{2\leq i\leq r}|V_{i}|^{2}\leq\frac{r-2}{r-1}(\sum_{2\leq i\leq r}|V_{i}|)^{2}=\frac{r-2}{r-1}(n-|V_{1}|)^{2}.$$
Then, noting $a=||V_{1}|-\frac{n}{r}|$,
\begin{equation}
\begin{aligned}
e(G)&\leq(\sum_{1\leq i<j\leq r}|V_{i}||V_{j}|)+(\sum_{1\leq i\leq r}e(G[V_{i}]))\\
&\leq|V_{1}|(n-|V_{1}|)+(\sum_{2\leq i<j\leq r}|V_{i}||V_{j}|)+\theta^{3}n^{2}\\
&\leq|V_{1}|(n-|V_{1}|)+\frac{r-2}{2(r-1)}(n-|V_{1}|)^{2}+\theta^{3}n^{2}\\
&=\frac{r-1}{2r}n^{2}-\frac{r}{2(r-1)}a^{2}+\theta^{3}n^{2}.
\end{aligned}\notag
\end{equation}
Recall that $e(G)\geq(\frac{r-1}{2r}-\theta^{3})n^{2}$. Thus, $(\frac{r-1}{2r}-\theta^{3})n^{2}\leq\frac{r-1}{2r}n^{2}-\frac{r}{2(r-1)}a^{2}+\theta^{3}n^{2}$, implying that $a\leq\sqrt{\frac{4(r-1)}{r}\theta^{3}n^{2}}\leq\theta n$. This finishes the proof. \hfill$\Box$

\medskip

For any $1\leq i\leq r$, let $W_{i}=\left\{v\in V_{i}~|~d_{V_{i}}(v)\geq\theta n\right\}$, and let $\overline{V}_{i}=V_{i}-W_{i}$.
Set $W=\cup_{1\leq i\leq r}W_{i}$.

\medskip

\f{\bf Claim 2.} $|W|\leq2\theta^{2} n$.

\medskip

\f{\bf Proof of Claim 2.} Since $\sum_{1\leq i\leq r}e(G[V_{i}])\leq\theta^{3}n^{2}$, and
$$\sum_{1\leq i\leq r}e(G[V_{i}])=\sum_{1\leq i\leq r}\frac{1}{2}\sum_{v\in V_{i}}d_{V_{i}}(v)\geq\sum_{1\leq i\leq r}\frac{1}{2}\sum_{v\in W_{i}}d_{V_{i}}(v)\geq\frac{1}{2}\sum_{1\leq i\leq r}\theta n|W_{i}|=\frac{1}{2}\theta n|W|,$$
we have $\frac{1}{2}\theta n|W|\leq\theta^{3}n^{2}$, implying that $|W|\leq2\theta^{2}n$.
 \hfill$\Box$

\medskip

\f{\bf Claim 3.}
Let $1\leq \ell\leq r$ be fixed. Assume that $u_{1},u_{2},...,u_{2rt^{2}}\in \cup_{1\leq i\neq \ell\leq r}\overline{V}_{i}$. For any constant $D_{0}>t$, there are $D_{0}$ vertices in $\overline{V}_{\ell}$ which are the common neighbors of $u_{1},u_{2},...,u_{2rt^{2}}$ in $G$.

\f{\bf Proof of Claim 3.} For any $1\leq i\leq 2rt^{2}$, assume that $u_{i}\in \overline{V}_{j_{i}}$, where $j_{i}\neq\ell$. Then $d_{G}(u_{i})\geq(\frac{r-1}{r}-\theta)n$ and $d_{V_{j_{i}}}(u_{i})\leq\theta n$ as $u_{i}\notin W$. Recall that $|V_{s}|\leq\frac{n}{r}+\theta n$ for any $1\leq s\leq r$.
Hence $$d_{V_{\ell}}(u_{i})= d_{G}(u_{i})-d_{V_{j_{i}}}(u_{i})-\sum_{1\leq s\neq\ell,j_{i}\leq r}d_{V_{s}}(u_{i})\geq(\frac{1}{r}-r\theta)n.$$
Then, using Fact 1, we have
$$|\cap_{1\leq i\leq 2rt^{2}}N_{V_{\ell}}(u_{i})|\geq(\sum_{1\leq i\leq 2rt^{2}}|N_{V_{\ell}}(u_{i})|)-(2rt^{2}-1)|V_{\ell}|\geq(\frac{1}{r}-2r(r+1)t^{2}\theta)n.$$
It follows that 
$$|\cap_{1\leq i\leq 2rt^{2}}N_{\overline{V}_{\ell}}(u_{i})|\geq(\frac{1}{r}-2r(r+1)t^{2}\theta)n-|W|
\geq(\frac{1}{r}-2r(r+2)t^{2}\theta)n>D_{0}.$$
Thus, there are $D_{0}$ vertices in $\overline{V}_{\ell}$ which are all adjacent to the vertices $u_{1},u_{2},...,u_{2rt^{2}}$.
 \hfill$\Box$

Since $G$ is $P_{t}\otimes T((r-1)t,r-1)$-free, by Claim 3, we see $G[\overline{V}_{i}]$ is $P_{t}$-free for any $1\leq i\leq r$. The following claim can be obtained using the method in \cite{S}. For completeness, we give a proof here.

\medskip

\f{\bf Claim 4.} $|W|\leq M$, where $M=\lceil4^{(r-1)t}(\frac{2}{\theta})^{t}rt\rceil$ is a constant. 

\f{\bf Proof of Claim 4.} Recall that $|W|\leq2\theta^{2}n$ by Claim 2.
For any fixed $1\leq i\leq r$ and $1\leq j\neq i\leq r$, we have $d_{V_{i}}(v)\leq d_{V_{j}}(v)$ for any $v\in V_{i}$, since $\sum_{1\leq \ell\leq r}e(G[V_{i}])$ is minimum (otherwise, we can obtain a contradiction by moving $v$ from $V_{i}$ to $V_{j}$). Therefore, $d_{V_{j}}(v)\geq\frac{1}{2}(d_{V_{i}}(v)+d_{V_{j}}(v))$ for any $v\in V_{i}$.  Recall that $|V_{s}|\leq(\frac{1}{r}+\theta)n$ for any $1\leq s\leq r$.
If $v\in W_{i}$, then
 $$d_{V_{j}}(v)\geq \frac{1}{2}(d_{G}(v)-\sum_{1\leq s\neq i,j\leq r}|V_{s}|)\geq \frac{1}{2}(\frac{r-1}{r}n-\theta n-(r-2)(\frac{1}{r}+\theta)n)\geq(\frac{1}{2r}-\frac{r-1}{2}\theta)n.$$
 It follows that
 $$d_{\overline{V}_{j}}(v)\geq d_{V_{j}}(v)-|W|\geq(\frac{1}{2r}-\frac{r-1}{2}\theta)n-2\theta^{2}n\geq\frac{n}{2r}-\frac{r}{2}\theta n.$$
If $v\in \overline{V}_{i}$, then $d_{V_{i}}(v)\leq \theta n$.
Thus
$$d_{V_{j}}(v)\geq d_{G}(v)-d_{V_{i}}(v)-\sum_{1\leq s\neq i,j\leq r}|V_{s}|\geq(\frac{r-1}{r}-\theta)n-\theta n-(r-2)(\frac{1}{r}+\theta)n\geq(\frac{1}{r}-r\theta)n.$$
It follows that
 $$d_{\overline{V}_{j}}(v)\geq d_{V_{j}}(v)-|W|\geq(\frac{1}{r}-r\theta)n-2\theta^{2}n\geq\frac{n}{r}-(r+1)\theta n.$$

In the following, we first show that $|W_{1}|\leq4^{(r-1)t}(\frac{2}{\theta})^{t}t$.

For any $v\in W_{1}$, we have that $d_{\overline{V}_{1}}(v)\geq d_{V_{1}}(v)-|W_{1}|\geq\theta n-2\theta^{2}n\geq\frac{1}{2}\theta n$. For $j\neq1$, we have $d_{\overline{V}_{j}}(v)\geq\frac{n}{2r}-\frac{r}{2}\theta n$ by the above discussion.
 Let $$Y_{1}=\left\{(w,Z)~|~w\in W_{1},Z\subseteq\overline{V}_{1},|Z|=t\right\},$$
  where $w$ is adjacent to the all the vertices in $Z$.
Since $d_{\overline{V}_{1}}(w)\geq\frac{1}{2}\theta n$ for any $w\in W_{1}$,  we have $|Y_{1}|\geq\tbinom{\frac{1}{2}\theta n}{t}|W_{1}|$. On the other hand,  $\overline{V}_{1}$ has at most $\tbinom{(\frac{1}{r}+\theta) n}{t}$ subsets $Z$ of size $t$. Hence, there is a subset $Z_{1}\subseteq \overline{V}_{1}$ of size $t$ such that the vertices in $Z_{1}$ have at least $\frac{\tbinom{\frac{1}{2}\theta n}{t}|W_{1}|}{\tbinom{(\frac{1}{r}+\theta) n}{t}}\geq(\frac{\theta}{2})^{t}|W_{1}|$ common neighbors  in $W_{1}$. Let $W^{1}\subseteq W_{1}$ be the set of common neighbors of the vertices in $Z_{1}$. Then $|W^{1}|\geq(\frac{\theta}{2})^{t}|W_{1}|$.

Suppose that for some $1\leq s\leq r-1$, we have obtained several sets $Z_{i}\subseteq \overline{V}_{i}$ of size $t$ for any $1\leq i\leq s$ and $W^{s}\subseteq W_{1}$ , such that any two vertices from distinct parts of $W^{s},Z_{1},...,Z_{s}$ are adjacent in $G$ and $|W^{s}|\geq(\frac{1}{4})^{(s-1)t}(\frac{\theta}{2})^{t}|W_{1}|$. Let $V'_{s+1}\subseteq \overline{V}_{s+1}$ be the set of common neighbors of the vertices in $\cup_{1\leq i\leq s}Z_{i}$. Since $(\frac{n}{r}-2\theta)n\leq|\overline{V}_{s+1}|\leq(\frac{n}{r}+2\theta)n$ and $d_{\overline{V}_{s+1}}(u)\geq\frac{n}{r}-(r+1)\theta n$ for any $u\in\cup_{1\leq i\leq s}Z_{i}$ by the above discussion,
we have that $d_{\overline{V}_{s+1}}(u)\geq|\overline{V}_{s+1}|-(r+3)\theta n$ for any $u\in\cup_{1\leq i\leq s}Z_{i}$. Using Fact 1, we have
$$|V'_{s+1}|=|\bigcap_{u\in\cup_{1\leq i\leq s}Z_{i}}N_{\overline{V}_{s+1}}(u)|\geq|\overline{V}_{s+1}|-st(r+3)\theta n\geq|\overline{V}_{s+1}|-rt(r+3)\theta n.$$
For any $w\in W_{1}$, note that $d_{\overline{V}_{s+1}}(w)\geq\frac{n}{2r}-\frac{r}{2}\theta n$ as above. Then $d_{V'_{s+1}}(w)\geq\frac{n}{2r}-rt(r+4)\theta n$.

Note that $|V'_{s+1}|\leq(\frac{1}{r}+\theta)n$.
Let $$Y_{s+1}=\left\{(w,Z)~|~w\in W^{s},Z\subseteq V'_{s+1},|Z|=t\right\},$$
 where $w$ is adjacent to all the vertices in $Z$.
Since $d_{V'_{s+1}}(w)\geq\frac{n}{2r}-rt(r+4)\theta n\geq\frac{n}{3r}$ for any $w\in W^{s}$,  we have $|Y_{s+1}|\geq|W^{s}|\tbinom{\frac{n}{3r}}{t}$. On the other hand,  $V'_{s+1}$ has at most $\tbinom{(\frac{1}{r}+\theta) n}{t}$ subsets $Z$ of size $t$. Hence, there is a subset $Z_{s+1}\subseteq V'_{s+1}$ of size $t$ such that the vertices in $Z_{s+1}$ have at least $\frac{|W^{s}|\tbinom{\frac{n}{3r}}{t}}{\tbinom{(\frac{1}{r}+\theta) n}{t}}\geq(\frac{1}{4})^{t}|W^{s}|$ common neighbors  in $W^{s}$. Let $W^{s+1}\subseteq W^{s}$ be the set of common neighbors of the vertices in $Z_{s+1}$. Then $$|W^{s+1}|\geq(\frac{1}{4})^{t}|W^{s}|\geq(\frac{1}{4})^{st}(\frac{\theta}{2})^{t}|W_{1}|.$$
Therefore, by induction on $s$, we can obtain several sets $Z_{i}\subseteq \overline{V_{i}}$ of size $t$ for any $1\leq i\leq r$  and $W^{r}\subseteq W_{1}$, such that any two vertices from distinct parts of $W^{r},Z_{1},...,Z_{r}$ are adjacent in $G$ and $|W^{r}|\geq(\frac{1}{4})^{(r-1)t}(\frac{\theta}{2})^{t}|W_{1}|$. So, if $|W_{1}|\geq4^{(r-1)t}(\frac{2}{\theta})^{t}t$, then $|W^{r}|\geq t$. It follows that the subgraph induced by $W^{r},Z_{1},...,Z_{r}$ contains a complete $(r+1)$-partite graph of which each part has size at least $t$. This contradicts the fact that $G$ is $\mathcal{F}$-free. Thus $|W_{1}|<4^{(r-1)t}(\frac{2}{\theta})^{t}t$. Similarly,  we can show that $|W_{i}|<4^{(r-1)t}(\frac{2}{\theta})^{t}t$ for any $2\leq i\leq r$. Hence $|W|<4^{(r-1)t}(\frac{2}{\theta})^{t}rt$.
This finishes the proof. \hfill$\Box$

\medskip

\medskip

For any fixed $1\leq i\leq r$, let $Q$ be an induced subgraph of $G[\overline{V}_{i}]$ (or a subset of $\overline{V}_{i}$).
Define
$$\alpha(Q)=E(Q)\cup\left\{uv\in E(G)~|~u\in V(Q),v\in (V_{i}-V(Q))\cup W\right\},$$
and $$\beta(Q)=\left\{uv\notin E(G)~|~u\in V(Q),v\in\cup_{1\leq j\neq i\leq r}\overline{V}_{j}\right\}.$$

\begin{lem}\label{symmetry edge bound}
For any fixed $1\leq i_{0}\leq r$, if $Q_{1},Q_{2},...,Q_{t}$ are symmetric subgraphs of $G[\overline{V}_{i_{0}}]$, then $|\alpha(Q_{i})|\leq(M+2t)|Q_{1}|$ and $|\beta(Q_{1})|\leq(M+2t)|Q_{1}|$. 
\end{lem}

\f{\bf Proof:} For convenience, assume $i_{0}=1$. We first show that there are at most $t$ vertices in $\overline{V}_{1}-\cup_{1\leq i\leq t}V(Q_{i})$ which have neighbors in $V(Q_{1})$. If not, there are $t$ vertices $u_{1},u_{2},...,u_{t}$ in $\overline{V}_{1}-\cup_{1\leq i\leq t}V(Q_{i})$, such that $u_{i}v_{i}$ is an edge for any $1\leq i\leq t$, where $v_{1},v_{2},...,v_{t}\in V(Q_{1})$ (may be repeated). Since $Q_{1},Q_{2},...,Q_{t}$ are symmetric subgraphs in $G[\overline{V}_{1}]$, for $2\leq j\leq t$ let $v^{j}_{1},v^{j}_{2},...,v^{j}_{t}\in V(Q_{j})$ be the symmetric vertices of $v_{1},v_{2},...,v_{t}$ in $V(Q_{j})$, respectively. Then $u_{i}$ is adjacent to $v^{j}_{i}$ for any $1 \leq i,j\leq t$. Note that $Q_{j}$ is connected for any $1\leq j\leq t$. Let $P_{v^{j}_{j}v^{j}_{j+1}}$ be a path connecting $v^{j}_{j}$ to $v^{j}_{j+1}$ in $Q_{j}$. Then
$$u_{1}P_{v^{1}_{1}v^{1}_{2}}u_{2}P_{v^{2}_{2}v^{2}_{3}}u_{3}\cdots u_{t-1}P_{v^{t-1}_{t-1}v^{t-1}_{t}}u_{t}$$
is a path of order at least $2t-1$ in $G[\overline{V}_{1}]$. Now we select its a sub-path of order $t$, say $w_{1},w_{2},...,w_{t}$. By Claim 3, these $t$ vertices have $t$ common neighbors in $\overline{V}_{2}$, say $w^{2}_{1},w^{2}_{2},...,w^{2}_{t}$. By Claim 3, these $2t$ vertices have $t$ common neighbors in $\overline{V}_{3}$. Repeat this process. Eventually, we can obtain a subgraph in $G$, which contains $P_{t}\otimes T(t(r-1),r-1)$ as a subgraph, a contradiction. Hence there are at most $t$ vertices in $\overline{V}_{1}-\cup_{1\leq i\leq t}V(Q_{i})$ which have neighbors in $V(Q_{1})$. Moreover, using a similar discussion, we see that $Q_{i}$ is $\left\{P_{t}\right\}$-free for any $1\leq i\leq t$. This implies that $e(Q_{1})\leq t|Q_{1}|$ by Lemma \ref{ex n Pt}. Hence, $|\alpha(Q_{1})|\leq (M+2t)|Q_{1}|$.

Suppose $|\beta(Q_{1})|>(M+2t)|Q_{1}|$. Let $G'$ be the graph obtained from $G$ by deleting all the edges in $\alpha(Q_{1})$, and adding all the non-edges in $\beta(Q_{1})$.
Then $F\subseteq G'$ for some $F\in\mathcal{F}$, since $e(G')>e(G)$ and $G\in {\rm EX}(n,\mathcal{F})$. Let $U$ be the set of vertices of $F$ which are adjacent to vertices in $V(F)\cap V(Q_{1})$ in $G'$. Then $U\subseteq(\cup_{2\leq i\leq r}\overline{V}_{i})$. By Claim 3, in $G$ we can choose $t$ common neighbors in $\overline{V}_{1}-V(F)$ of the vertices in $U$, say $y_{1},y_{2},...,y_{t}$. Let $F'$ be the graph of $G$ induced by $(V(F)-V(F)\cap V(Q_{1}))\cup\left\{y_{1},y_{2},...,y_{t}\right\}$. Clearly, $F\subseteq F'$, a contradiction as $G$ is $F$-free.
Hence, $|\beta(Q_{1})|\leq(M+2t)|Q_{1}|$. For any $Q_{i}$ with $2\leq i\leq t$, we can use a similar discussion.
This finishes the proof.
 \hfill$\Box$

\medskip

The idea of Lemma \ref{symmetry extended} should be attributed to Simonovits \cite{S1}. 

\begin{lem}\label{symmetry extended}
 Let $1\leq i_{0}\leq r$ be fixed. Assume that $Q_{1},Q_{2},...,Q{\tau}$ are symmetric subgraphs of $G[V_{i_{0}}\cup (W\cup B)]$, where $B\subseteq V(G)-(V_{i_{0}}\cup W)$ and $V(Q_{i})\subseteq\overline{V}_{i_{0}}$ for any $1\leq i\leq\tau$. If $Q_{1},Q_{2},...,Q{\tau}$ are not symmetric in $G$, then there exists one vertex $v\notin V_{i_{0}}\cup (W\cup B)$, such  that at least $\frac{1}{2^{|Q_{1}|+1}t}\tau$ ones among the $\tau$ subgraphs, say $Q_{1},Q_{2},...,Q{\frac{1}{2^{|Q_{1}|+1}t}\tau}$,  are also symmetric in $G[V_{i_{0}}\cup (W\cup B\cup\left\{v\right\})]$, and $v$ has at least one non-neighbor in $Q_{i}$ for any $1\leq i\leq\frac{1}{2^{|Q_{1}|+1}t}\tau$.
\end{lem}

\f{\bf Proof:}  We can assume that $\frac{1}{2^{|Q_{1}|+1}t}\tau\geq1$. Otherwise, there is nothing to prove. Then $\tau\geq2^{|Q_{1}|+1}t$. For convenience, assume $i_{0}=1$. 
Let $Z=V(G)-(V_{1}\cup(W\cup B))$. If all vertices of $Q_{i}$ are adjacent to all the vertices in $Z$ for any $1\leq i\leq\tau$, then $Q_{1},Q_{2},...,Q_{\tau}$ are symmetric in $G$, too. This contradicts the assumption. Thus, without loss of generality, we can assume that $w_{1}\in Z$ has a non-neighbor $u_{1}\in V(Q_{1})$. Let $G'$ be the graph obtained from $G$ by adding the edge $u_{1}w_{1}$. Since $G\in{\rm EX}(n,\mathcal{F})$ and $e(G')>e(G)$, we see that $G'$ contains a subgraph $F\in \mathcal{F}$. Clearly, $F$ contains the edge $u_{1}w_{1}$.

Let $v_{1},...,v_{j}$ be all the vertices in $Z$ which are  contained in $F$, where $1\leq j\leq t$. Now we show that there are at most $t$ ones among $Q_{1},Q_{2},...,Q_{\tau}$, whose vertices are all adjacent to all the vertices $v_{1},...,v_{j}$.
Suppose not.  Noting $|F|\leq t$, we can find a $Q_{\ell}$ with $1\leq \ell\leq \tau$ and $V(Q_{\ell})\cap V(F)=\emptyset$ , such that all vertices of $Q_{\ell}$ are  adjacent to all the vertices $v_{1},...,v_{j}$. Recall that $Q_{1},Q_{2},...,Q_{\tau}$ are symmetric in $G[V_{1}\cup(W\cup B)]$. We see that all the neighbors in $F$ of the vertices in $F\cap V(Q_{1})$  are adjacent to the corresponding vertices in $V(Q_{\ell})$ in $G$. Let $F'$ be the subgraph of $G$ induced by $(V(F)-V(Q_{1}))\cup V(Q_{\ell})$. Clearly, $F\subseteq F'$. Thus $G$ contains a copy of $F$, a contradiction. Therefore, there are at most $t$ ones among $Q_{1},Q_{2},...,Q_{\tau}$, whose vertices are all adjacent to all the vertices $v_{1},...,v_{j}$. Hence, there are at least $\tau-t$ ones among $Q_{1},Q_{2},...,Q_{\tau}$, which have a non-neighbor in $\left\{v_{1},...,v_{j}\right\}$. So, we can find a vertex among $v_{1},...,v_{j}$, say $v$, such that $v$ has a non-neighbor in at least $\frac{\tau-t}{j}\geq\frac{\tau-t}{t}\geq\frac{1}{2t}\tau$ ones among $Q_{1},Q_{2},...,Q_{\tau}$, say $Q_{1},Q_{2},...,Q_{\frac{1}{2t}\tau}$. Clearly, $v$ is connecting in the same way to at least $\frac{1}{2t}\frac{1}{2^{|Q_{1}|}}\tau=\frac{1}{2^{|Q_{1}|+1}t}\tau$ ones among these $\frac{1}{2t}\tau$ subgraphs, say  $Q_{1},Q_{2},...,Q_{\frac{1}{2^{|Q_{1}|+1}t}\tau}$. Clearly, $Q_{1},Q_{2},...,Q_{\frac{1}{2^{|Q_{1}|+1}t}\tau}$ are also symmetric in $G[V_{i_{0}}\cup (W\cup B\cup\left\{v\right\})]$.
This finishes the proof.
 \hfill$\Box$

\begin{lem}\label{symmetry to G}
 Let $1\leq i_{0}\leq r$ be fixed. Assume that $Q_{1},Q_{2},...,Q_{\tau}$ are symmetric subgraphs of  $G[V_{i_{0}}\cup W]$, where $V(Q_{i})\subseteq\overline{V}_{i_{0}}$ for any $1\leq i\leq\tau$. Then there are at least $(\frac{1}{2^{|Q_{1}|+1}t})^{(M+2t)|Q_{1}|}\cdot\tau$ ones among $Q_{1},Q_{2},...,Q_{\tau}$ are symmetric in $G$, too.
\end{lem}

\f{\bf Proof:} We can assume that $(\frac{1}{2^{|Q_{1}|+1}t})^{(M+2t)|Q_{1}|}\tau\geq1$. Otherwise, there is nothing to prove. Then $\tau\geq(2^{|Q_{1}|+1}t)^{(M+2t)|Q_{1}|}$. For convenience, assume $i_{0}=1$. By Lemma \ref{symmetry edge bound}, $|\beta(Q_{i})|\leq(M+2t)|Q_{1}|$ for any $1\leq i\leq\tau$.  

 If $Q_{1},Q_{2},...,Q_{\tau}$ are symmetric in $G$, then the lemma holds. So, we can assume that $Q_{1},Q_{2},...,Q_{\tau}$ are not symmetric in $G$. By Lemma \ref{symmetry extended}, there exists a vertex $v_{1}\notin V_{1}\cup W$, such that at least $\frac{1}{2^{|Q_{1}|+1}t}\tau$ ones among the $\tau$ subgraphs, say $Q_{1},Q_{2},...,Q_{\frac{1}{2^{|Q_{1}|+1}t}\tau}$, which are symmetric in $G[V_{1}\cup (W\cup\left\{v_{1}\right\})]$. Moreover, $v_{1}$ has a non-neighbor in $V(Q_{i})$ for any $1\leq i\leq\frac{1}{2^{|Q_{1}|+1}t}\tau$. 

 If $Q_{1},Q_{2},...,Q_{\frac{1}{2^{|Q_{1}|+1}t}\tau}$ are symmetric in $G$, then the lemma holds. So, we can assume that $Q_{1},Q_{2},...,Q_{\frac{1}{2^{|Q_{1}|+1}t}\tau}$ are not symmetric in $G$. By Lemma \ref{symmetry extended} again, there exists another vertex $v_{2}\notin V_{1}\cup(W\cup \left\{v_{1}\right\})$, such that at least $(\frac{1}{2^{|Q_{1}|+1}t})^{2}\tau$ ones among the subgraphs, say $Q_{1},Q_{2},...,Q_{(\frac{1}{2^{|Q_{1}|+1}t})^{2}\tau}$, which are symmetric in $G[V_{1}\cup (W\cup\left\{v_{1},v_{2}\right\})]$. Moreover, $v_{2}$ has a non-neighbor in $V(Q_{i})$ for any $1\leq i\leq(\frac{1}{2^{|Q_{1}|+1}t})^{2}\tau$. 

Repeat the above process. We can obtain vertices $v_{1}, v_{2},..., v_{\ell}\notin V_{1}\cup W$, such that $Q_{1},Q_{2},...,Q_{(\frac{1}{2^{|Q_{1}|+1}t})^{\ell}\cdot\tau}$ are symmetric in $G[V_{1}\cup(W\cup \left\{v_{1}, v_{2},..., v_{\ell}\right\})]$, and each vertex in $\left\{v_{1}, v_{2},..., v_{\ell}\right\}$ has a non-neighbor in $V(Q_{i})$ for any $1\leq i\leq(\frac{1}{2^{|Q_{1}|+1}t})^{\ell}\cdot\tau$. Since $|\beta(Q_{i})|\leq(M+2t)|Q_{1}|$ for any $1\leq i\leq\tau$, we see that $\ell\leq(M+2t)|Q_{1}|$. Since the process must stop at step $\ell\leq(M+2t)|Q_{1}|$, we must have that $Q_{1},Q_{2},...,Q_{(\frac{1}{2^{|Q_{1}|+1}t})^{\ell}\cdot\tau}$ are symmetric in $G$, too. The proof is completed as 
$$(\frac{1}{2^{|Q_{1}|+1}t})^{\ell}\cdot\tau\geq
(\frac{1}{2^{|Q_{1}|+1}t})^{(M+2t)|Q_{1}|}\cdot\tau.$$
 \hfill$\Box$

\medskip

Now we come back to the proof Lemma \ref{Pl-free}.
Let $f(t)=\sum_{0\leq i\leq t-1}t^{i}$. Clearly, a connected graph $Q$ must contain $K_{1,t}$ or $P_{t}$ as a subgraph if $|Q|\geq f(t)$. 
Set $$N=1\times2\times\cdots\times f(t)\times N_{0}.$$

 Let $\pi=\left\{m_{k}\right\}_{k\geq1}$, where $m_{k}=2^{(f(t)(M+2t)r+M)k}\cdot(2^{k+1}t)^{(M+2t)k}\cdot N$ for any $k\geq1$. Set $D=D(t,\pi)+f(t)$. Let $\pi'=\left\{m'_{k}\right\}_{k\geq1}$, where $m'_{k}=2^{(D(M+2t)r+M)k}\cdot(2^{k+1}t)^{(M+2t)k}\cdot D\cdot N$ for any $k\geq1$. By Corollary \ref{symmetric}, if a $P_{t}$-free graph $H$ has order at least $C(t,\pi)$, then $H$ has $2^{(f(t)(M+2t)r+M)k}\cdot(2^{k+1}t)^{(M+2t)k}\cdot N$ symmetric subgraphs of order $k$, where $1\leq k\leq D$. 
If a $P_{t}$-free graph $H$ has order at least $C(t,\pi')$, then $H$ has $2^{(D(M+2t)r+M)k}\cdot(2^{k+1}t)^{(M+2t)k}\cdot D\cdot N$ symmetric subgraphs of order $k$, where $1\leq k\leq D(t,\pi')$. Set $C=\max\left\{C(t,\pi),C(t,\pi')\right\}$.
In our proof, we can assume that $n\gg C$.

\medskip

\f{\bf Case 1.} For some $1\leq i_{0}\leq r$, there are $N$ symmetric subgraphs $Q_{1},Q_{2},...,Q_{N}$ of $G$, such that  $\cup_{1\leq j\leq N}V(Q_{j})\subseteq \overline{V}_{i_{0}}$ and $f(t)\leq |Q_{1}|\leq D$.

\medskip

Recall that $Q_{1},Q_{2},...,Q_{N}$ are symmetric in $G$ with $f(t)\leq|Q_{1}|\leq D$. Without loss of generality, assume that $i_{0}=1$. That is, $V(Q_{i})\subseteq \overline{V}_{1}$ for any $1\leq i\leq N$.  By Corollary \ref{symmetry edge bound}, $|\beta(Q_{i})|\leq(M+2t)|Q_{i}|\leq(M+2t)D$ for any $1\leq i\leq N$. Let $\overline{N}(Q_{1})$ be the set of vertices in $\overline{V}_{2}$ which are not adjacent to at least one vertex in $V(Q_{1})$. Then $|\overline{N}(Q_{1})|\leq|\beta(Q_{1})|\leq(M+2t)D$. Let $\overline{V}'_{2}=\overline{V}_{2}-\overline{N}(Q_{1})$. Then by symmetry, all the vertices in $\cup_{1\leq i\leq N}V(Q_{i})$ are adjacent to all the vertices in $\overline{V}'_{2}$.
 
Since $|\overline{V}'_{2}|\geq (\frac{1}{r}-\theta)n-M-(M+2t)D>C$, by Corollary \ref{symmetric}, $G[\overline{V}'_{2}]$ has  $2^{(D(M+2t)r+M)k}\cdot(2^{k+1}t)^{(M+2t)k}\cdot D\cdot N$ symmetric subgraphs $Q^{2}_{1},Q^{2}_{2},Q^{2}_{3},...$,  where $|Q^{2}_{1}|=k\leq D(t,\pi')$. Set $\tau=2^{(D(M+2t)r+M)|Q^{2}_{1}|}\cdot(2^{|Q^{2}_{1}|+1}t)^{(M+2t)|Q^{2}_{1}|}\cdot D\cdot N$. Now we prove that $|Q^{2}_{1}|< f(t)$. Suppose not, i.e., $|Q^{2}_{1}|\geq f(t)$. Then $Q^{2}_{j}$ contains a $K_{1,t}$ for any $1\leq j\leq t$,  since $Q^{2}_{j}$ is $P_{t}$-free. Recall that $|Q_{j}|\geq f(t)$. thus, $Q_{j}$ contains a $K_{1,t}$ for any $1\leq j\leq t$. Hence, we obtain a subgraph $A_{1}=tK_{1,t}$ in  $G[\overline{V}_{1}]$ and a subgraph $A_{2}=tK_{1,t}$ in $G[\overline{V}_{2}]$. Moreover, the vertices of  $A_{1}$ are all adjacent to all the vertices of  $A_{2}$. By Claim 3, the vertices of $A_{1},A_{2}$ have $t$ common neighbors in $\overline{V}_{3}$, say $u^{3}_{1},u^{3}_{2},...,u^{3}_{t}$.  By Claim 3 again, the vertices of $A_{1},A_{2},\left\{u^{3}_{1},u^{3}_{2},...,u^{3}_{t}\right\}$ have $t$ common neighbors in $\overline{V}_{4}$, say $u^{4}_{1},u^{4}_{2},...,u^{4}_{t}$.  Repeat the above process. We can obtain a subgraph $tK_{1,t}\otimes tK_{1,t}\otimes T(t(r-2),r-2)$ in $G$. Hence $|Q^{2}_{1}|\leq f(t)$. For any $1\leq j\leq\tau$, there are at most $2^{|Q^{2}_{1}|((M+2t)D+M)}$ ways connecting the vertices of $Q^{2}_{j}$ to the vertices of $W\cup \overline{N}(Q_{1})$. Thus, among $Q^{2}_{j}$ with $1\leq j\leq\tau$, there are at least $2^{-|Q^{2}_{1}|((M+2t)D+M)}\cdot\tau$ ones which are also symmetric in $G[V_{2}\cup W]$. Denote these ones by $Q^{2}_{1},Q^{2}_{2},...,Q^{2}_{2^{-|Q^{2}_{1}|((M+2t)D+M)}\cdot\tau}$ without loss of generality. By Lemma \ref{symmetry to G}, among $Q^{2}_{j}$ with $1\leq j\leq2^{-|Q^{2}_{1}|((M+2t)D+M)}\cdot\tau$, there are at least $(\frac{1}{2^{|Q^{2}_{1}|+1}t})^{(M+2t)|Q^{2}_{1}|}\cdot2^{-|Q^{2}_{1}|((M+2t)D+M)}\cdot\tau\geq D\cdot N$ ones which are also symmetric in $G$. Note that $N$ is divisible by $|Q^{2}_{1}|$ as $|Q^{2}_{1}|\leq f(t)$. It follows that $\frac{|Q_{1}|\cdot N}{|Q^{2}_{1}|}\leq D\cdot N$ is an integer. Thus, without loss of generality, we can assume that $Q^{2}_{1},Q^{2}_{2},...,Q^{2}_{\frac{|Q_{1}|\cdot N}{|Q^{2}_{1}|}}$ are symmetric in $G$.

 Let $\overline{N}(Q_{1}\cup Q^{2}_{1})$ be the set of vertices in $\overline{V}_{3}$ which are not adjacent to at least one vertex in $V(Q_{1}\cup Q^{2}_{1})$. Then $|\overline{N}(Q_{1}\cup Q^{2}_{1})|\leq|\beta(Q_{1}\cup Q^{2}_{1})|\leq2(M+2t)D$. Let $\overline{V}'_{3}=\overline{V}_{3}-\overline{N}(Q_{1}\cup Q^{2}_{1})$. Then by symmetry, all the vertices in $(\cup_{1\leq i\leq N}V(Q_{i}))\cup(\cup_{1\leq i\leq \frac{|Q_{1}|\cdot N}{|Q^{2}_{1}|}}V(Q^{2}_{i}))$ are adjacent to all the vertices in $\overline{V}'_{3}$.
 Since $|\overline{V}'_{3}|\geq (\frac{1}{r}-\theta)n-M-2(M+2t)D>C$, by Corollary \ref{symmetric}, $G[\overline{V}'_{3}]$ has  $2^{(D(M+2t)r+M)|Q^{3}_{1}|}\cdot(2^{|Q^{3}_{1}|+1}t)^{(M+2t)|Q^{3}_{1}|}\cdot D\cdot N$ symmetric subgraphs $Q^{3}_{1},Q^{3}_{2},Q^{3}_{3},...$,  where $|Q^{3}_{1}|\leq D(t,\pi')$. Set $\tau_{1}=2^{(D(M+2t)r+M)|Q^{3}_{1}|}\cdot(2^{|Q^{3}_{1}|+1}t)^{(M+2t)|Q^{3}_{1}|}\cdot D\cdot N$. Similar to the previous paragraph, we can prove that $|Q^{3}_{1}|< f(t)$. Moreover, among $Q^{3}_{j}$ with $1\leq j\leq\tau_{1}$, we can assume that $Q^{3}_{1},Q^{3}_{2},...,Q^{3}_{\frac{|Q_{1}|\cdot N}{|Q^{3}_{1}|}}$ are symmetric in $G$.

Repeat the above process. We can find symmetric subgraphs in $G$: 
$$Q_{1},Q_{2},...,Q_{N},$$
$$Q^{2}_{1},Q^{2}_{2},...,Q^{2}_{\frac{|Q_{1}|\cdot N}{|Q^{2}_{1}|}},$$
$$Q^{3}_{1},Q^{3}_{2},...,Q^{3}_{\frac{|Q_{1}|\cdot N}{|Q^{3}_{1}|}},$$
$$\cdots\cdots$$
$$Q^{r}_{1},Q^{r}_{2},...,Q^{r}_{\frac{|Q_{1}|\cdot N}{|Q^{r}_{1}|}}.$$
Moreover, all the vertices from distinct groups of symmetric subgraphs are adjacent in $G$. 
Each group of symmetric subgraphs has exactly $|Q_{1}|N$ vertices. Let $|Q_{1}|N=O(r,t)\cdot N_{0}$. Clearly, $O(r,t)$ is an integer, and $|Q_{1}|,|Q^{j}_{1}|\leq O(r,t)$ for any $2\leq j\leq r$, as desired.

\medskip

\f{\bf Case 2.} For any $1\leq i\leq r$, if $Q_{1},Q_{2},...,Q_{N}$ are symmetric subgraphs of $G$ such that $\cup_{1\leq j\leq N}V(Q_{j})\subseteq \overline{V}_{i}$, then $|Q_{1}|< f(t)$ or $|Q_{1}|>D$.

\medskip

Ignore the symbols used in Case 1. Since $|\overline{V}_{1}|\geq (\frac{1}{r}-\theta)n-M>C$, by Corollary \ref{symmetric}, $G[\overline{V}_{1}]$ has  $2^{(f(t)(M+2t)r+M)k}\cdot(2^{k+1}t)^{(M+2t)k}\cdot N$ symmetric subgraphs $Q^{1}_{1},Q^{1}_{2},Q^{1}_{3},...$,  where $|Q^{1}_{1}|=k\leq D$. Set $\tau_{1}=2^{(f(t)(M+2t)r+M)|Q^{1}_{1}|}\cdot(2^{|Q^{1}_{1}|+1}t)^{(M+2t)|Q^{1}_{1}|}\cdot N$. 
 For any $1\leq j\leq\tau_{1}$, there are at most $2^{|Q^{1}_{1}|M}$ ways connecting the vertices of $Q^{1}_{j}$ to the vertices of $W$. Thus, among $Q^{1}_{j}$ with $1\leq j\leq\tau_{1}$, there are at least $2^{-|Q^{1}_{1}|M}\cdot\tau_{1}$ ones which are also symmetric in $G[V_{1}\cup W]$. Denote these ones by $Q^{1}_{1},Q^{1}_{2},...,Q^{1}_{2^{-|Q^{1}_{1}|M}\cdot\tau_{1}}$ without loss of generality. By Lemma \ref{symmetry to G}, among $Q^{1}_{j}$ with $1\leq j\leq2^{-|Q^{1}_{1}|M}\cdot\tau_{1}$, there are at least $(\frac{1}{2^{|Q^{1}_{1}|+1}t})^{(M+2t)|Q^{1}_{1}|}\cdot2^{-|Q^{1}_{1}|M}\cdot\tau_{1}\geq N$ ones which are also symmetric in $G$. Denote these ones by $Q^{1}_{1},Q^{1}_{2},...,Q^{1}_{N}$ without loss of generality. Since $|Q^{1}_{1}|\leq D$, we must have $|Q^{1}_{1}|\leq f(t)$ by the assumption in this case. 

Let $\overline{N}(Q^{1}_{1})$ be the set of vertices in $\overline{V}_{2}$ which are not adjacent to at least one vertex in $V(Q^{1}_{1})$. Then $|\overline{N}(Q^{1}_{1})|\leq|\beta(Q^{1}_{1})|\leq(M+2t)|Q^{1}_{1}|\leq f(t)(M+2t)$. Let $\overline{V}'_{2}=\overline{V}_{2}-\overline{N}(Q^{1}_{1})$. Then by symmetry, all the vertices in $\cup_{1\leq i\leq N}V(Q^{1}_{i})$ are adjacent to all the vertices in $\overline{V}'_{2}$.
Since $|\overline{V}'_{2}|\geq (\frac{1}{r}-\theta)n-M-f(t)(M+2t)>C$, by Corollary \ref{symmetric},  $G[\overline{V}'_{2}]$ has  $2^{(f(t)(M+2t)r+M)k}\cdot(2^{k+1}t)^{(M+2t)k}\cdot N$ symmetric subgraphs $Q^{2}_{1},Q^{2}_{2},Q^{2}_{3},...$,  where $|Q^{2}_{1}|=k\leq D$. Set $\tau_{2}=2^{(f(t)(M+2t)r+M)|Q^{2}_{1}|}\cdot(2^{|Q^{2}_{1}|+1}t)^{(M+2t)|Q^{2}_{1}|}\cdot N$. 
 For any $1\leq j\leq\tau_{2}$, there are at most $2^{|Q^{2}_{1}|(f(t)(M+2t)+M)}$ ways connecting the vertices of $Q^{2}_{j}$ to the vertices of $W\cup\overline{N}(Q^{1}_{1})$. Thus, among $Q^{1}_{j}$ with $1\leq j\leq\tau_{2}$, there are at least $2^{-|Q^{2}_{1}|(f(t)(M+2t)+M)}\cdot\tau_{2}$ ones which are also symmetric in $G[V_{2}\cup W]$. Denote these ones by $Q^{2}_{1},Q^{2}_{2},...,Q^{2}_{2^{-|Q^{2}_{1}|(f(t)(M+2t)+M)}\cdot\tau_{2}}$ without loss of generality. By Lemma \ref{symmetry to G}, among $Q^{1}_{j}$ with $1\leq j\leq2^{-|Q^{2}_{1}|(f(t)(M+2t)+M)}\cdot\tau_{2}$, there are at least $(\frac{1}{2^{|Q^{1}_{1}|+1}t})^{(M+2t)|Q^{1}_{1}|}
 \cdot2^{-|Q^{2}_{1}|(f(t)(M+2t)+M)}\cdot\tau_{2}\geq N$ ones which are also symmetric in $G$. Denote these ones by $Q^{2}_{1},Q^{2}_{2},...,Q^{2}_{N}$ without loss of generality. Since $|Q^{2}_{1}|\leq D$, we must have $|Q^{2}_{1}|\leq f(t)$ by the assumption in this case.

 Let $\overline{N}(Q^{1}_{1}\cup Q^{2}_{1})$ be the set of vertices in $\overline{V}_{3}$ which are not adjacent to at least one vertex in $V(Q^{2}_{1}\cup Q^{2}_{1})$. Then $|\overline{N}(Q^{1}_{1}\cup Q^{2}_{1})|\leq|\beta(Q^{1}_{1}\cup Q^{2}_{1})|\leq2(M+2t)f(t)$. Let $\overline{V}'_{3}=\overline{V}_{3}-\overline{N}(Q^{1}_{1}\cup Q^{2}_{1})$. Then by symmetry, all the vertices in $(\cup_{1\leq i\leq N}V(Q^{1}_{i}))\cup(\cup_{1\leq i\leq N}V(Q^{2}_{i}))$ are adjacent to all the vertices in $\overline{V}'_{3}$.
 Since $|\overline{V}'_{3}|\geq (\frac{1}{r}-\theta)n-M-2(M+2t)f(t)>C$, by Corollary \ref{symmetric}, $G[\overline{V}'_{3}]$ has  $2^{(f(t)(M+2t)r+M)|Q^{3}_{1}|}\cdot(2^{|Q^{3}_{1}|+1}t)^{(M+2t)|Q^{3}_{1}|}\cdot N$ symmetric subgraphs $Q^{3}_{1},Q^{3}_{2},Q^{3}_{3},...$,  where $|Q^{3}_{1}|\leq D$. Similar to the previous paragraph, we can find symmetric subgraphs of $G$: $Q^{3}_{1},Q^{3}_{2},...,Q^{3}_{N}$, such that $|Q^{3}_{1}|\leq f(t)$. 

Repeat the above process. We can find symmetric subgraphs in $G$: 
$$Q^{1}_{1},Q^{1}_{2},...,Q^{1}_{N},$$
$$Q^{2}_{1},Q^{2}_{2},...,Q^{2}_{N},$$
$$\cdots\cdots$$
$$Q^{r}_{1},Q^{r}_{2},...,Q^{r}_{N},$$
where $|Q^{j}_{1}|\leq f(t)$ for any $1\leq j\leq r$.
Moreover, all the vertices from distinct groups of symmetric subgraphs are adjacent in $G$. 
Without loss of generality, assume that $|Q_{1}|=\min_{1\leq j\leq r}|Q_{j}|$. Note that $N$ is divisible by $|Q_{j}|$. Thus, $\frac{|Q_{1}|N}{|Q_{j}|}\leq N$ is an integer for any $1\leq j\leq r$. For any $2\leq j\leq r$, we only need $Q^{j}_{1},Q^{j}_{2},...,Q^{j}_{\frac{|Q_{1}|N}{|Q_{j}|}}$. Hence, we can obtain $r$ new groups of symmetric subgraphs of $G$, in which each group has exactly $|Q_{1}|N$ vertices. Let $|Q_{1}|N=O(r,t)\cdot N_{0}$. Clearly, $O(r,t)$ is an integer, and $|Q^{j}_{1}|\leq O(r,t)$ for any $1\leq j\leq r$, as desired.
This completes the proof. \hfill$\Box$

\medskip

\section{Proof of Theorem \ref{matching-free}}

In this section, we give the proof of Theorem \ref{matching-free}.

\f{\bf Proof of Theorem \ref{matching-free}.} Recall $q=q(\mathcal{F})$. For large $m$, Let $\mathcal{G}_{m}$ be the set of $\mathcal{F}$-free graphs $H$ such that, $H$ has a partition $V(H)=W\cup(\cup_{1\leq i\leq r}S_{i})$ satisfying:\\
$(i)$ $|W|=q-1$ and $|S_{i}|=\lfloor\frac{n-q+1}{r}\rfloor$ or $\lceil\frac{n-q+1}{r}\rceil$ for any $1\leq i\leq r$;\\
$(ii)$ for any $1\leq i\leq r$, there exists $S'_{i}\subseteq S_{i}$ such that $|S_{i}-S'_{i}|\leq t^{2}$, and $N_{H}(v)=V(G)-S_{i}$ for any $v\in S'_{i}$.

Note that $\mathcal{G}_{m}$ is not empty, since $\overline{K_{q-1}}\otimes T(n-q+1,r)\in\mathcal{G}_{m}$ by the definition of $q$.

For the graphs in $\mathcal{G}_{m}$, we can define an operation $\mathcal{D}$ as follows: 
Let $H$ be a graph in $\mathcal{G}_{m}$ with the partition as above. Define $\mathcal{D}(H)$ to be the graph obtained from $H$ by adding a new vertex $x_{i}$ to $S_{i}$ for each $1\leq i\leq r$, such that $N_{\mathcal{D}(H)}(x_{i})=V(\mathcal{D}(H))-(S_{i}\cup\left\{x_{i}\right\})$. Clearly, $\mathcal{D}(H)\in\mathcal{G}_{m+r}$, and $e(\mathcal{D}(H))=e(H)+(q-1)+m(r-1)+\binom{r}{2}$.
The inverse of $\mathcal{D}$ is defined as follows:
Let $H$ be a graph in $\mathcal{G}_{m}$ with the partition as above. Define $\mathcal{D}^{-1}(H)$ to be the graph obtained from $H$ by deleting a vertex in $S_{i}$ for each $1\leq i\leq r$. Clearly, $\mathcal{D}^{-1}(H)\in\mathcal{G}_{m-r}$, and $e(\mathcal{D}^{-1}(H))=e(H)-(q-1)+(m-r)(r-1)+\binom{r}{2}$.
Note that $\mathcal{D}^{-1}(\mathcal{D}(H))=\mathcal{D}(\mathcal{D}^{-1}(H))=H$ for any $H\in\mathcal{G}_{m}$. For any $\ell\geq1$, define iteratively $\mathcal{D}^{\ell+1}(H)=\mathcal{D}(\mathcal{D}^{\ell}(H))$, where $\mathcal{D}^{1}(H)=\mathcal{D}(H)$.

\medskip

\f{\bf Claim 1.} Let $H$ be an extremal graph with maximum edges in $\mathcal{G}_{m}$. Then $\mathcal{D}^{\pm1}(H)$ is an extremal graph with maximum edges in $\mathcal{G}_{m\pm r}$.

\medskip

\f{\bf Proof of Claim 1.} For similarity, we only prove the result for $\mathcal{D}^{-1}(H)$. Suppose that $\mathcal{D}^{-1}(H)$ is not an extremal graph in $\mathcal{G}_{m-r}$. Then there is a graph $H'$ in $\mathcal{G}_{m-r}$ with $e(H')>e(\mathcal{D}^{-1}(H))$. Then $$\mathcal{D}(H')=e(H')+(q-1)+(m-r)(r-1)+\binom{r}{2}$$
$$>e(\mathcal{D}^{-1}(H))+(q-1)+(m-r)(r-1)+\binom{r}{2}=e(H).$$
Note that $H,\mathcal{D}(H')\in\mathcal{G}_{m}$. This contradicts the fact that $H$ is extremal in $\mathcal{G}_{m}$. Hence, $\mathcal{D}^{-1}(H)$ is extremal in $\mathcal{G}_{m-r}$.
This finishes the proof of Claim 1. \hfill$\Box$

\medskip

For any constant $\ell\geq1$ and any extremal graph $H$ in $\mathcal{G}_{m}$ (with large $m$), by Claim 1, $\mathcal{D}^{-\ell}(H)$ is extremal in $\mathcal{G}_{m-r\ell}$. We shall prove that any graph in ${\rm EX}(n,\mathcal{F})$ with large $n$ is in $\mathcal{G}_{n}$ using the progressive induction. 

For $m\geq1$,  let $\mathfrak{U}_{m}={\rm EX}(m,\mathcal{F})$.  For a graph $G'\in{\rm EX}(m,\mathcal{F})$, we say $G'$ has property $B$, if $G'\in \mathcal{G}_{m}$. Choose an $H_{m}$ as an extremal graph in $\mathcal{G}_{m}$. For any $G_{m}\in{\rm EX}(n,\mathcal{F})$,
define $\phi(G_{m})=e(G_{m})-e(H_{m})$. Clearly, $\phi(G_{m})$ is a non-negative integer-valued function. If $G_{m}$ has property $B$, then $\phi(G_{m})=0$. By Lemma \ref{progressive induction}, it suffices to show that for sufficiently large $m$, either $G_{m}$ satisfies property $B$, or there is a $\frac{m}{2}\leq m'<m$, such that $G_{m'}\in{\rm EX}(m',\mathcal{F})$ and $\phi(G_{m'})>\phi(G_{m})$.

Recall that $F_{1}\subseteq tK_{2}\otimes T(t(r-1),r-1)$,  where $F_{1}\in\mathcal{F}$. Since $|F_{1}|\leq t$, we see $F_{1}\subseteq P_{t}\otimes T(t(r-1),r-1)$ and $F_{1}\subseteq tK_{1,t}\otimes tK_{1,t}\otimes T(t(r-2),r-2)$. Thus $\mathcal{F}$ satisfies the condition in Lemma \ref{Pl-free}. Since $G_{m}\in {\rm EX}(m,\mathcal{F})$ with large $m$, there is a constant $N$ (by letting $N_{0}=2t$ and  $N=O(r,t)\cdot N_{0}$ in Lemma \ref{Pl-free}), such that $G_{m}$  has an induced subgraph $\otimes_{1\leq i\leq r}G^{i}$ satisfying:\\
 $(1)$  $|G^{i}|=N$ for any $1\leq i\leq r$;\\
 $(2)$ for any $1\leq i\leq r$, $G^{i}=\cup _{1\leq j\leq \ell_{i}}H^{i}_{j}$, where $H^{i}_{1},H^{i}_{2},...,H^{i}_{\ell_{i}}$ are symmetric subgraphs of $G$ and $\ell_{i}\geq 2t$.
 
Since $G_{m}$ is $F_{1}$-free, we must have that $|H^{i}_{1}|=1$ for any $1\leq i\leq r$. That is, $G^{i}$ consisting of symmetric vertices of $G$ for any $1\leq i\leq r$. 
 
Let $V'=V(G)-V(\otimes_{1\leq i\leq r}G^{i})$. By symmetry, for each $G^{i}$, a vertex $u\in V'$ is  adjacent either to all the vertices of $G^{i}$, or to no vertex of $G^{i}$. Now we divide the vertices in $V'$ by several parts. Let $X$ be the set of vertices $u$ in $V'$, such that $u$ is adjacent to all the vertices in $\cup_{1 \leq i\leq r}V(G^{i})$. For $1\leq i\leq r$, let $X_{i}$ be the set of vertices $u$ in $V'$, such that $u$ is adjacent to all the vertices in $\cup_{1 \leq j\neq i\leq r}V(G^{j})$, but to no vertex of $V(G^{i})$. Let $Y=V'-(X\cup(\cup_{1\leq i\leq r}X_{i}))$. Note that each vertex in $Y$ is adjacent to at most $(r-2)N$ vertices of $\otimes_{1\leq i\leq r}G^{i}$. By the definition of $q$, we must have $|X|\leq q-1$. By a calculation. we have $e(H_{m})=e(\mathcal{D}^{-N}(H_{m}))+e(T(rN,r))+(q-1)N+(m-rN)(r-1)N$. Note that $e(\mathcal{D}^{-N}(H_{m}))=e(H_{m-rN})$ by Claim 1. Thus,
\begin{equation}
\begin{aligned}
\phi(G_{m})&=e(G_{m})-e(H_{m})\\
&=e(G_{m}[V'])+e_{G}(V(\otimes_{1\leq i\leq r}G^{i}),V')+e(T(rN,r))\\
&-(e(H_{m-rN})+e(T(rN,r))+(q-1)N+(m-rN)(r-1)N)\\
&\leq\phi(G_{m-rN})+e_{G}(V(\otimes_{1\leq i\leq r}G^{i}),V')-((q-1)N+(m-rN)(r-1)N)).
\end{aligned}\notag
\end{equation}
If $e_{G}(V(\otimes_{1\leq i\leq r}G^{i}),V')-((q-1)N+(m-rN)(r-1)N))<0$, then $\phi(G_{m})<\phi(G_{m-rN}),$ as desired (noting $m-rN\geq\frac{m}{2}$ for large $m$).
 Thus, we can assume that 
 $$e_{G}(V(\otimes_{1\leq i\leq r}G^{i}),V')\geq((q-1)N+(m-rN)(r-1)N)).$$
  We will show that $G_{m}$ has property $B$.
Note that
 \begin{equation}
\begin{aligned}
&e_{G}(V(\otimes_{1\leq i\leq r}G^{i}),V')\\
&=e_{G}(V(\otimes_{1\leq i\leq r}G^{i}),X)+e_{G}(V(\otimes_{1\leq i\leq r}G^{i}),\cup_{1\leq i\leq r}X_{i})+e_{G}(V(\otimes_{1\leq i\leq r}G^{i}),Y)\\
&\leq |X|rN+(r-1)N|\cup_{1\leq i\leq r}X_{i}|+(r-2)N|Y|\\
&=|X|N+(r-1)N(m-rN)-N|Y|\\
&\leq (q-1)N+(r-1)N(m-rN)-N|Y|.
\end{aligned}\notag
\end{equation}
 Thus, we must have $|X|=q-1$ and $Y=\emptyset$.
 
 For each $1\leq i\leq r$, let $S_{i}=V(G^{i})\cup X_{i}$. Then $G_{m}$ has a partition $V(G_{m})=X\cup(\cup_{1\leq i\leq r}S_{i})$ satisfying:\\
$(3)$ $|X|=q-1$ and $S_{i}=V(G^{i})\cup X_{i}$ for any $1\leq i\leq r$, where $|V(G^{i})|=N\geq2t$;\\
$(4)$ for any $1\leq i\leq r$, $N_{G_{m}}(v)=V(G)-S_{i}$ for any $v\in V(G^{i})$.

Recall that $F_{1}\subseteq tK_{2}\otimes T(t(r-1),r-1),F_{2}\subseteq (K_{q-1,t}\cup K_{1,t})\otimes T(t(r-1),r-1)$, and $G_{m}$ is $F_{1}$-free and $F_{2}$-free. Now we show that, for any $1\leq i\leq r$, $G_{m}[X_{i}]$ is $tK_{2}$-free and $K_{1,t}$-free. For similarity, we only show this result for $i=1$. If $G_{m}[X_{1}]$ contains a $tK_{2}$, then $G_{m}[X_{1}\cup(\cup_{2\leq j\leq r}V(G^{j}))]$ contains a $F_{1}$, a contradiction. If $G_{m}[X_{1}]$ contains a $K_{1,t}$, then $G_{m}[(X\cup S_{1})\cup(\cup_{2\leq j\leq r}V(G^{j}))]$ contains a $F_{2}$, a contradiction. Hence, $G_{m}[X_{1}]$ is $tK_{2}$-free and $K_{1,t}$-free.
 It follows that $G_{m}[X_{1}]$ will contain no edges after deleting at most $t^{2}$ vertices.
Let $S'_{1}$ be the subset of $S_{1}$ after deleting these $t^{2}$ vertices. Then $|S_{1}-S'_{1}|\leq t^{2}$. For any $v\in S'_{1}-X_{1}$, if $v$ is not adjacent to some vertex $w\in V(G_{m})-S_{1}$, let $G''$ be the graph obtained from $G_{m}$ by adding all the non-edges between $v$ and $ V(G_{m})-S_{1}$. Note that $v$ and the vertices of $V(G^{1})$ are symmetric in $G''$. By Observation 1, we see that $G''$ is $\mathcal{F}$-free. But $e(G'')>e(G_{m})$, a contradiction as $G_{m}\in{\rm EX}(m,\mathcal{F})$. Thus $N_{G_{m}}(v)=V(G_{m})-S_{1}$ for any $v\in S'_{1}$. For $2\leq i\leq r$, $S'_{i}$ is similar to obtain. 

Now we show that $||S_{i}|-|S_{j}||\leq1$ for any $1\leq i<j\leq r$. Suppose not. Without loss of generality, assume that $|S_{1}|\geq|S_{2}|+2$. let $G'''$ be the graph obtained from $G_{m}$ by moving a vertex $v\in S'_{1}$ to $S_{2}$, i.e., deleting all the edges incident with $v$ in $G_{m}$ and connecting $v$ to all the vertices in $ V(G_{m})-S_{2}$. Note that $v$ and the vertices of $S'_{2}$ are symmetric in $G'''$. Again by Observation 1, we see that $G'''$ is $\mathcal{F}$-free. But $e(G''')>e(G_{m})$, a contradiction as $G_{m}\in{\rm EX}(m,\mathcal{F})$. Thus $||S_{i}|-|S_{j}||\leq1$ for any $1\leq i<j\leq r$. It follows that $|S_{i}|=\lfloor\frac{n-q+1}{r}\rfloor$ or $\lceil\frac{n-q+1}{r}\rceil$ for any $1\leq i\leq r$.
This completes the proof.\hfill$\Box$

\section{Proofs of Theorem \ref{general Icosahe} and Theorem \ref{m Icosahe}}

In this section, we give the proofs of Theorem \ref{general Icosahe} and Theorem \ref{m Icosahe}. We first prove Theorem \ref{general Icosahe}.

\f{\bf Proof of Theorem \ref{general Icosahe}.} 
$(i).$ Recall that $2\leq r\leq\ell-3$ ($\ell$ is even), and $\mathcal{F}=\left\{F_{1},F_{2},F_{3}\right\}$, where
$$F_{1}=(P_{\ell}\cup\overline{K_{m}})\otimes T(m(r-1),r-1),$$
$$F_{2}=(K_{1,2}\cup mK_{2})\otimes mK_{2}\otimes T(m(r-2),r-2),$$ $$F_{3}=(K_{\frac{\ell}{2}}\cup K_{\ell-1}\cup\overline{K_{m}})\otimes T(m(r-1),r-1).$$
Clearly, $\min_{1\leq i\leq3}\chi(F_{i})=r+1$. Moreover, for any $1\leq i\leq 3$, $\chi(F_{i}-U)\geq r+1$ for any $U\subseteq F_{i}$ with $|U|=\frac{\ell-2}{2}$. It follows that 
$K_{\frac{\ell-2}{2}}\otimes T(n-\frac{\ell-2}{2},r)$ is $\mathcal{F}$-free.
We will prove that ${\rm EX}(n,\mathcal{F})=\left\{K_{\frac{\ell-2}{2}}\otimes T(n-\frac{\ell-2}{2},r)\right\}$ for sufficiently large $n$ using the progressive induction.

For $n\geq1$,  let $\mathfrak{U}_{n}={\rm EX}(n,\mathcal{F})$.
For a  graph $G'\in{\rm EX}(n,\mathcal{F})$, we say $G'$ has property $B$, if $G'=K_{\frac{\ell-2}{2}}\otimes T(n-\frac{\ell-2}{2},r)$.  For a  graph $G\in{\rm EX}(n,\mathcal{F})$, define 
$$\phi(G)=e(G)-e(K_{\frac{\ell-2}{2}}\otimes T(n-\frac{\ell-2}{2},r)).$$
 Clearly, $\phi(G)$ is a non-negative integer-valued function. If $G$ has property $B$, then $\phi(G)=0$. By Lemma \ref{progressive induction}, it suffices to show that for sufficiently large $n$, either $G$ satisfies property $B$, or there is a $\frac{n}{2}\leq n'<n$ and a $G'\in{\rm EX}(n',\mathcal{F})$ such that $\phi(G')>\phi(G)$.

Since $F_{1},F_{2}\in\mathcal{F},$ we see that $\mathcal{F}$ satisfies the condition in Lemma \ref{Pl-free} (noting $t=\max_{1\leq i\leq3}|F_{i}|$ there). Since $G\in {\rm EX}(n,\mathcal{F})$ with sufficiently large $n$, there is a constant $N$ (by letting $N_{0}=3m(\ell-1)$ and  $N=O(r,t)\cdot N_{0}$ in Lemma \ref{Pl-free}), such that $G$  has an induced subgraph $\otimes_{1\leq i\leq r}G^{i}$ satisfying:\\
 $(1)$  $|G^{i}|=N$ for any $1\leq i\leq r$;\\
 $(2)$ for any $1\leq i\leq r$, $G^{i}=\cup _{1\leq j\leq \ell_{i}}H^{i}_{j}$, where $H^{i}_{1},H^{i}_{2},...,H^{i}_{\ell_{i}}$ are symmetric subgraphs of $G$ and $\ell_{i}\geq 3m$.\\
 Note that $3m(\ell-1)| N$.

 \medskip
 
 \f{\bf Case 1.} There are at least two $i$, such that $|H^{i}_{1}|\geq2$.
 
\medskip

In this case, we can assume that $|H^{i}_{1}|\geq2$ for $i=1,2$ without loss of generality. Note that $\ell_{i}\geq3m$. Since $G$ is $F_{2}$-free, we must have that $|H^{1}_{1}|=|H^{2}_{1}|=2$ and $|H^{i}_{1}|\leq2$ for any $3\leq i\leq r$. That is, $G^{i}=\ell_{i}K_{1}$ or $G^{i}=\ell_{i}K_{2}$ for any $1\leq i\leq r$. 
 
Let $V'=V(G)-V(\otimes_{1\leq i\leq r}G^{i})$. By symmetry, for $G^{i}=\ell_{i}K_{1}$, a vertex $u\in V'$ is  adjacent either to all the vertices of $G^{i}$, or to no vertex of $G^{i}$. For $G^{i}=\ell_{i}K_{2}$, a vertex $u\in V'$ is  adjacent either to all the vertices of $G^{i}$, or to half number of the vertices of $G^{i}$, or to no vertex of $G^{i}$. Suppose that there is a vertex $v\in V'$ such that $v$ is adjacent to at least $(r-1)N+1$ vertices of $\otimes_{1\leq i\leq r}G^{i}$. Then by symmetry, $v$ must be adjacent to all the vertices of $\otimes_{1\leq i\neq i_{0}\leq r}G^{i}$ and to at least half number of the vertices of $G^{i_{0}}$ ($=\ell_{i_{0}}K_{2}$) for some $1\leq i_{0}\leq r$. Without loss of generality, assume $i_{0}=1$. Then $G[\left\{v\right\}\cup V(G^{1}\otimes G^{2})]$ contains a $(K_{1,2}\cup mK_{2})\otimes mK_{2}$. And then $G[\left\{v\right\}\cup V(\otimes_{1\leq i\leq r}G^{i})]$ contains a $F_{2}$, a contradiction.
Thus, for any $v\in V'$, $v$ is adjacent to at most $(r-1)N$ vertices of $\otimes_{1\leq i\leq r}G^{i}$. It follows that $e_{G}(V(\otimes_{1\leq i\leq r}G^{i}),V')\leq(n-rN)(r-1)N$.

By a calculation. we have 
$e(K_{\frac{\ell-2}{2}}\otimes T(n-\frac{\ell-2}{2},r))=e(K_{\frac{\ell-2}{2}}\otimes T(n-rN-\frac{\ell-2}{2},r))+e(T(rN,r))+\frac{\ell-2}{2}N+(n-rN)(r-1)N.$
 Let $G'\in{\rm EX}(n-rN,\mathcal{F})$. Then $e(G[V'])\leq e(G')$. Note that $\phi(G')=e(G')-e(K_{\frac{\ell-2}{2}}\otimes T(n-rN-\frac{\ell-2}{2},r))$.
 Thus,
\begin{equation}
\begin{aligned}
\phi(G)&=e(G)-e(K_{\frac{\ell-2}{2}}\otimes T(n-\frac{\ell-2}{2},r))\\
&\leq e(G[V'])+e_{G}(V(\otimes_{1\leq i\leq r}G^{i}),V')+e(T(rN,r))+\frac{N}{2}r\\
&-e(K_{\frac{\ell-2}{2}}\otimes T(n-rN-\frac{\ell-2}{2},r))-e(T(rN,r))-\frac{\ell-2}{2}N-(n-rN)(r-1)N\\
&\leq\phi(G')+(n-rN)(r-1)N+\frac{N}{2}r-\frac{\ell-2}{2}N-(n-rN)(r-1)N\\
&=\phi(G')-\frac{N}{2}(\ell-2-r)~(noting~r\leq\ell-3)\\
&<\phi(G'),
\end{aligned}\notag
\end{equation}
as desired.

\medskip
 
 \f{\bf Case 2.} There is at most one $i$, such that $|H^{i}_{1}|\geq2$.
 
\medskip

In this case, we can assume that $|H^{i}_{1}|=1$ for $i\geq2$ without loss of generality. That is $G_{i}=N\cdot K_{1}$ for $i\geq2$. Recall that $3m(\ell-1)|N$.
 Also, let $V'=V(G)-V(\otimes_{1\leq i\leq r}G^{i})$. By symmetry, we can divide the vertices in $V'$ by $r+1$ parts as follows. Let $X_{1}$ be the set of vertices $u$ in $V'$, such that $u$ is adjacent to all the vertices of $\otimes_{2 \leq i\leq r}G^{i}$. For $2\leq i\leq r$, let $X_{i}$ be the set of vertices $u$ in $V'$, such that $u$ is adjacent to all the vertices of $\otimes_{1 \leq j\neq i\leq r}G^{j}$ but to no vertex of $V(G^{i})$. Let $Y=V'-\cup_{1\leq i\leq r}X_{i}$, i.e., the set of vertices which have non-neighbors in at least two among $G^{i}$ with $1\leq i\leq r$. Note that each vertex in $\cup_{2\leq i\leq r}X_{i}$ is adjacent to  $(r-1)N$ vertices of $\otimes_{1\leq i\leq r}G^{i}$. Each vertex in $Y$ is adjacent to at most $(r-1)N-1$ vertices of $\otimes_{1\leq i\leq r}G^{i}$ (using symmetry). The way connecting the vertices in $X_{1}$ to $V(G^{1})$ is not known now. 
 
Recall that $e(K_{\frac{\ell-2}{2}}\otimes T(n-\frac{\ell-2}{2},r))=e(K_{\frac{\ell-2}{2}}\otimes T(n-rN-\frac{\ell-2}{2},r))+e(T(rN,r))+\frac{\ell-2}{2}N+(n-rN)(r-1)N$. Let $G'\in{\rm EX}(n-rN,\mathcal{F})$. Then $e(G[V'])\leq e(G')$. Note that $\phi(G')=e(G')-e(K_{\frac{\ell-2}{2}}\otimes T(n-rN-\frac{\ell-2}{2},r))$. 
Thus,
 \begin{equation}
\begin{aligned}
\phi(G)&=e(G)-e(K_{\frac{\ell-2}{2}}\otimes T(n-\frac{\ell-2}{2},r))\\
&=e(G[V'])+e_{G}(V(\otimes_{1\leq i\leq r}G^{i}),V')+e(T(rN,r))+e(G^{1})\\
&-e(K_{\frac{\ell-2}{2}}\otimes T(n-rN-\frac{\ell-2}{2},r))-e(T(rN,r))-\frac{\ell-2}{2}N-(n-rN)(r-1)N\\
&\leq\phi(G')+e_{G}(V(\otimes_{1\leq i\leq r}G^{i}),V')+e(G^{1})-\frac{\ell-2}{2}N-(n-rN)(r-1)N.
\end{aligned}\notag
\end{equation}
If $e_{G}(V(\otimes_{1\leq i\leq r}G^{i}),V')+e(G^{1})-\frac{\ell-2}{2}N-(n-rN)(r-1)N<0$, then $\phi(G)<\phi(G')$, as desired. 
 Thus, we can assume that 
 $$e_{G}(V(\otimes_{1\leq i\leq r}G^{i}),V')+e(G^{1})\geq\frac{\ell-2}{2}N+(n-rN)(r-1)N.$$
  We will show that $G$ has property $B$ (i.e., $G=K_{\frac{\ell-2}{2}}\otimes T(n-\frac{\ell-2}{2},r)$).

Note that
 \begin{equation}
\begin{aligned}
&e_{G}(V(\otimes_{1\leq i\leq r}G^{i}),V')\\
&=e_{G}(V(\otimes_{1\leq i\leq r}G^{i}),X_{1})+e_{G}(V(\otimes_{1\leq i\leq r}G^{i}),\cup_{2\leq i\leq r}X_{i})+e_{G}(V(\otimes_{1\leq i\leq r}G^{i}),Y)\\
&\leq e_{G}(V(\otimes_{1\leq i\leq r}G^{i}),X_{1})+(r-1)N|\cup_{2\leq i\leq r}X_{i}|+((r-1)N-1)|Y|\\
&=e_{G}(V(\otimes_{1\leq i\leq r}G^{i}),X_{1})+(r-1)N\cdot(n-rN-|X_{1}|)-|Y|\\
&=e_{G}(V(G^{1}),X_{1})+(r-1)N|X_{1}|+(r-1)N\cdot(n-rN-|X_{1}|)-|Y|\\
&=e_{G}(V(G^{1}),X_{1})+(r-1)N\cdot(n-rN)-|Y|.
\end{aligned}\notag
\end{equation}
Recall $e_{G}(V(\otimes_{1\leq i\leq r}G^{i}),V')+e(G^{1})\geq\frac{\ell-2}{2}N+(n-rN)(r-1)N$.
It follows that 
$$e(G^{1})+e_{G}(V(G^{1}),X_{1})\geq\frac{\ell-2}{2}N+|Y|\geq\frac{\ell-2}{2}N.$$
Note that each vertex in $V(G^{1})\cup X_{1}$ is adjacent to all the vertices in $V(\otimes_{2\leq i\leq r}G^{i})$. Thus, $G[V(G^{1})\cup X_{1}]$ is $P_{\ell}$-free, since $G$ is $F_{1}$-free. Note that $H^{1}_{1},H^{1}_{1},...,H^{1}_{\ell_{1}}$ are also symmetric in $G[V(G^{1})\cup X_{1}]$ and $|G^{1}|=N$. By Lemma \ref{desired symmetric subgraphs}, we have
$e(G^{1})+e_{G}(V(G^{1}),X_{1})=\frac{\ell-2}{2}N$ (implying $|Y|=0$).
Moreover, either $Q^{1}_{i}=K_{\ell-1}$ for any $1\leq i\leq\ell_{1}$ and $e_{G}(V(G^{1}),X_{1})=0$, or $Q^{1}_{i}=K_{1}$ for any $1\leq i\leq\ell_{1}$ and there are $\frac{\ell-2}{2}$ vertices in $X_{1}$ which are all adjacent to the (unique) vertex in $Q^{1}_{i}$ for all $1\leq i\leq\ell_{1}$. If $Q^{1}_{i}=K_{\ell-1}$ for any $1\leq i\leq\ell_{1}$, then $F_{3}$ is contained in $G$, a contradiction. Hence, we must have the latter case. That is, $G^{1}=N\cdot K_{1}$, and there are $\frac{\ell-2}{2}$ vertices in $X_{1}$, say $y_{1},y_{2},...,y_{\frac{\ell-2}{2}}$, such that they are all adjacent to all the vertices of $G^{1}$. Recall $Y=\emptyset$. Let $X=\left\{y_{1},y_{2},...,y_{\frac{\ell-2}{2}}\right\}$. For convenience, we still use $X_{1}$ to denote $X_{1}-X$.
 For each $1\leq i\leq r$, let $S_{i}=V(G^{i})\cup X_{i}$. Then $G$ has a partition $V(G)=X\cup(\cup_{1\leq i\leq r}S_{i})$ satisfying:\\
$(3)$ $X=\left\{y_{1},y_{2},...,y_{\frac{\ell-2}{2}}\right\}$ and $S_{i}=V(G^{i})\cup X_{i}$ for any $1\leq i\leq r$, where $|V(G^{i})|=N\geq3m$;\\
$(4)$ for any $1\leq i\leq r$, $N_{G}(v)=V(G)-S_{i}$ for any $v\in V(G^{i})$.

For any $1\leq i\leq r$, let $X'_{i}$ be the set of vertices in $S_{i}$ which have neighbors in $X$. Now we show that the vertices in $X'_{i}$ are isolated in $G[S_{i}]$. Suppose not. Without loss of generality, assume that $v\in X'_{1}$ such that $v$ is adjacent to $y_{\frac{\ell-2}{2}}$ and adjacent to another vertex $w$ in $S_{1}$. Let $u_{1},u_{2},...,u_{\frac{\ell-2}{2}}$ be other $\frac{\ell-2}{2}$ vertices (not $v,w$) in $V(G^{1})$. Then $u_{1}y_{1}u_{2}y_{2}\cdots u_{\frac{\ell-2}{2}}y_{\frac{\ell-2}{2}}vw$ is a $P_{\ell}$. Since the vertices of this path are all adjacent to the vertices in $V(\otimes_{2\leq i\leq r}G^{i})$, we see $F_{1}$ is contained in $G$, a contradiction. Hence, the vertices in $X'_{i}$ are isolated in $G[S_{i}]$ for any $1\leq i\leq r$.

Note that $G[X_{i}-X'_{i}]$ is $P_{\ell}$-free for any $1\leq i\leq r$. By Lemma \ref{ex n Pt}, we have $e(G[X_{i}-X'_{i}])\leq\frac{\ell-2}{2}|X_{i}-X'_{i}|$ with equality if and only if $G[X_{i}-X'_{i}]$ is the disjoint union of several copies of $K_{\ell-1}$. 
Thus,
\begin{equation}
\begin{aligned}
e(G)&\leq e(K_{\frac{\ell-2}{2}}\otimes K_{|S_{1}|,|S_{2}|,...,|S_{r}|})+(\sum_{1\leq i\leq r}e(G[X_{i}-X'_{i}]))-(\sum_{1\leq i\leq r}\frac{\ell-2}{2}|X_{i}-X'_{i}|)\\
&\leq e(K_{\frac{\ell-2}{2}}\otimes T(n-\frac{\ell-2}{2},r))+(\sum_{1\leq i\leq r}e(G[X_{i}-X'_{i}]))-(\sum_{1\leq i\leq r}\frac{\ell-2}{2}|X_{i}-X'_{i}|)\\
&\leq e(K_{\frac{\ell-2}{2}}\otimes T(n-\frac{\ell-2}{2},r)).
\end{aligned}\notag
\end{equation}
However, $e(G)\geq e(K_{\frac{\ell-2}{2}}\otimes T(n-\frac{\ell-2}{2},r))$ as $G\in{\rm EX}(n,\mathcal{F})$. This forces that $G[X]=K_{\frac{\ell-2}{2}}$, and $||S_{i}|-\frac{n-\frac{\ell-2}{2}}{r}|<1$ and $e(G[X_{i}-X'_{i}])=\frac{\ell-2}{2}|X_{i}-X'_{i}|$ for any $1\leq i\leq r$.
Then $G[X_{i}-X'_{i}]$ is the disjoint union of $K_{\ell-1}$. If $X_{i}-X'_{i}$ is not empty for some $1\leq i\leq r$, say $i=1$ without loss of generality, then $G[X_{1}-X'_{1}]$ contains a $K_{\ell-1}$. Now we choose a vertex $u_{1}$ in $V(G^{1})$. Then $G[\left\{u_{1},y_{1},y_{2},...,y_{\frac{\ell-2}{2}}\right\}]=K_{\frac{\ell}{2}}$. Clearly, $G[X\cup S_{1}]$ contains a $K_{\frac{\ell}{2}}\cup K_{\ell-1}\cup\overline{K_{m}}$.
Since all vertices of $\otimes_{2\leq i\leq r}G^{i}$ are adjacent to vertices in $X\cup S_{1}$, $G$ contains a $F_{3}$, a contradiction. Hence, $X_{i}-X'_{i}=\emptyset$ for any $1\leq i\leq r$. This means that all the vertices in $S_{i}$ are isolated in $G[S_{i}]$ for any $1\leq i\leq r$. Then, $G$ is a spanning subgraph of $K_{\frac{\ell-2}{2}}\otimes T(n-\frac{\ell-2}{2},r)$. Hence, $G=K_{\frac{\ell-2}{2}}\otimes T(n-\frac{\ell-2}{2},r)$, as desired.   

\medskip

$(ii)$. Ignore the symbols used in $(i)$. Recall that $2\leq r=\ell-2$ is even, and $\mathcal{F}=\left\{F_{1},F_{2},F_{3}\right\}$, where
$$F_{1}=(P_{r+2}\cup\overline{K_{m}})\otimes T(m(r-1),r-1),$$
$$F_{2}=(K_{1,2}\cup mK_{2})\otimes mK_{2}\otimes T(m(r-2),r-2),$$ $$F_{3}=(K_{\frac{r+2}{2}}\cup K_{r+1}\cup\overline{K_{m}})\otimes T(m(r-1),r-1).$$
Clearly, $\min_{1\leq i\leq3}\chi(F_{i})=r+1$. Moreover, for any $1\leq i\leq 3$, $\chi(F_{i}-U)\geq r+1$ for any $U\subseteq F_{i}$ with $|U|=\frac{r}{2}$. It follows that 
$K_{\frac{r}{2}}\otimes T(n-\frac{r}{2},r)$ is $\mathcal{F}$-free. For $n\geq2r$, recall that $G_{n,r}$
  is the graph obtained from $T(n,r)$ by adding a (near) perfect matching in each of its $r$ parts. By a result in \cite{S}, we see that $G_{n,r}$ is $K_{1}\otimes(\otimes_{1\leq i\leq r}2K_{1})$-free. Clearly, $K_{1}\otimes(\otimes_{1\leq i\leq r}2K_{1})$ is contained in $F_{i}$ for any $1\leq i\leq3$. Thus, $G_{n,r}$ is $\mathcal{F}$-free. We will prove that ${\rm EX}(n,\mathcal{F})\subseteq\left\{K_{\frac{r}{2}}\otimes T(n-\frac{r}{2},r)\right\}\cup\mathcal{G}_{n,r}$ for sufficiently large $n$ using the progressive induction.

For $n\geq1$,  let $\mathfrak{U}_{n}={\rm EX}(n,\mathcal{F})$.
For a  graph $G'\in{\rm EX}(n,\mathcal{F})$, we say $G'$ has property $B$, if $G'\in \left\{K_{\frac{r}{2}}\otimes T(n-\frac{r}{2},r)\right\}\cup\mathcal{G}_{n,r}$.  For a  graph $G\in{\rm EX}(n,\mathcal{F})$, define 
$$\phi(G)=\min\left\{e(G)-e(K_{\frac{r}{2}}\otimes T(n-\frac{r}{2},r)),e(G)-e(G_{n,r})\right\}.$$
 Clearly, $\phi(G)$ is a non-negative integer-valued function. If $G$ has property $B$, then $\phi(G)=0$. By Lemma \ref{progressive induction}, it suffices to show that for sufficiently large $n$, either $G$ satisfies property $B$, or there is a $\frac{n}{2}\leq n'<n$ and a $G'\in{\rm EX}(n',\mathcal{F})$ such that $\phi(G')>\phi(G)$.

 For even integer $N$, it is clear that (for large $n$) $e(G_{n,r})=e(G_{n-2N,r})+e(T(rN,r))+\frac{rN}{2}+(n-rN)(r-1)N$, and
  $e(K_{\frac{r}{2}}\otimes T(n-\frac{r}{2},r))=e(K_{\frac{r}{2}}\otimes T(n-rN-\frac{r}{2},r))+e(T(rN,r))+\frac{r}{2}N+(n-rN)(r-1)N$. For $G\in{\rm EX}(n,\mathcal{F})$ and $G'\in{\rm EX}(n-rN,\mathcal{F})$,
  it follows that
\begin{equation}
\begin{aligned}
\phi(G)&=\min\left\{e(G)-e(K_{\frac{r}{2}}\otimes T(n-\frac{r}{2},r)),e(G)-e(G_{n,r})\right\}\\
&=e(G)-e(G')+\min\left\{e(G')-e(K_{\frac{r}{2}}\otimes T(n-\frac{r}{2},r)),e(G')-e(G_{n,r})\right\}\\
&=e(G)-e(G')-(e(T(rN,r))+\frac{r}{2}N+(n-rN)(r-1)N)\\
&+\min\left\{e(G')-e(K_{\frac{r}{2}}\otimes T(n-rN-\frac{r}{2},r)),e(G')-e(G_{n-rN,r})\right\}\\
&=\phi(G')+e(G)-e(G')-(e(T(rN,r))+\frac{r}{2}N+(n-rN)(r-1)N).
\end{aligned}\notag
\end{equation}
If $e(G)-e(G')-(e(T(rN,r))+\frac{r}{2}N+(n-rN)(r-1)N)<0$, then $\phi(G)<\phi(G')$, as desired.  Thus, we can assume 
$$e(G)\geq e(G')+e(T(rN,r))+\frac{r}{2}N+(n-rN)(r-1)N$$
 in the following. We will show that $G$ has property $B$ (i.e., $G\in \left\{K_{\frac{r}{2}}\otimes T(n-\frac{r}{2},r)\right\}\cup\mathcal{G}_{n,r}$).
Recall that $G'\in{\rm EX}(n-rN,\mathcal{F})$.

Since $F_{1},F_{2}\in\mathcal{F}$, we see that $\mathcal{F}$ satisfies the condition in Lemma \ref{Pl-free}. Since $G\in {\rm EX}(n,\mathcal{F})$ with sufficiently large $n$, there is a constant $N$ (by letting $N_{0}=4m$ and  $N=O(r,t)\cdot N_{0}$ in Lemma \ref{Pl-free}), such that $G$  has an induced subgraph $\otimes_{1\leq i\leq r}G^{i}$ satisfying:\\
 $(1)$  $|G^{i}|=N$ for any $1\leq i\leq r$;\\
 $(2)$ for any $1\leq i\leq r$, $G^{i}=\cup _{1\leq j\leq \ell_{i}}H^{i}_{j}$, where $H^{i}_{1},H^{i}_{2},...,H^{i}_{\ell_{i}}$ are symmetric subgraphs of $G$ and $\ell_{i}\geq 4m$.\\
 Note that $4m| N$.

 \medskip
 
 \f{\bf Case 1.} There are at least two $i$, such that $|H^{i}_{1}|\geq2$.
 
\medskip

In this case, we can assume that $|H^{i}_{1}|\geq2$ for $i=1,2$ without loss of generality. Note that $\ell_{i}\geq4m$. Since $G$ is $F_{2}$-free, we must have that $|H^{1}_{1}|=|H^{2}_{1}|=2$ and $|H^{i}_{1}|\leq2$ for any $3\leq i\leq r$. 
Let $V'=V(G)-V(\otimes_{1\leq i\leq r}G^{i})$. Suppose that there is a vertex $v\in V'$ such that $v$ is adjacent to at least $(r-1)N+1$ vertices of $\otimes_{1\leq i\leq r}G^{i}$. Then by symmetry, $v$ must be adjacent to all the vertices of $\otimes_{1\leq i\neq i_{0}\leq r}G^{i}$ and to at least half number of the vertices of $G^{i_{0}}$ for some $1\leq i_{0}\leq 2$. Clearly, $G[\left\{v\right\}\cup V(\otimes_{1\leq i\leq r}G^{i})]$ contains a $F_{2}$, a contradiction.
Thus, for any $v\in V'$, $v$ is adjacent to at most $(r-1)N$ vertices of $\otimes_{1\leq i\leq r}G^{i}$. Let $Y$ be the set of vertices $v\in V'$, such that $v$ is adjacent to at most $(r-1)N-1$ vertices of $\otimes_{1\leq i\leq r}G^{i}$. It follows that $e_{G}(V(\otimes_{1\leq i\leq r}G^{i}),V')\leq(n-rN)(r-1)N-|Y|$.

Since $G[V']$ is $\mathcal{F}$-free, we have $e(G[V'])\leq e(G')$ as $G'\in{\rm EX}(n-rN,\mathcal{F})$. Recall that $e(G)\geq e(G')+e(T(rN,r))+\frac{r}{2}N+(n-rN)(r-1)N$. Thus,
$e(G)\geq e(G[V'])+e(T(rN,r))+\frac{r}{2}N+(n-rN)(r-1)N$. Note that $e(\otimes_{1\leq i\leq r}G^{i})\leq\frac{rN}{2}+e(T(rN,r))$.
 However,
\begin{equation}
\begin{aligned}
e(G)&=e(G[V'])+e_{G}(V(\otimes_{1\leq i\leq r}G^{i}),V')+e(\otimes_{1\leq i\leq r}G^{i})\\
&\leq e(G[V'])+(n-rN)(r-1)N-|Y|+\frac{rN}{2}+e(T(rN,r)).
\end{aligned}\notag
\end{equation}
We must have $|Y|=0,e(\otimes_{1\leq i\leq r}G^{i})=\frac{rN}{2}+e(T(rN,r))$, and each vertex in $V'$ is adjacent to exactly $(r-1)N$ vertices of $\otimes_{1\leq i\leq r}G^{i}$. It follows that $G^{i}=\frac{N}{2}K_{2}$ for any $1\leq i\leq r$. If some vertex $v\in V'$ has neighbors in $V(G^{i})$ for each $1\leq i\leq r$, then by symmetry, $G$ contains $F_{1}$ for $r=2$, and $G$ contains $F_{2}$ for $r\geq4$ ($r$ is even), a contradiction. Hence, each vertex $v\in V'$ is adjacent to no vertex of $V(G^{i})$ for some $1\leq i\leq r$, and thus adjacent to all vertices of $\otimes_{1\leq j\neq i\leq r}G^{j}$. Then by the above discussion, $V'=\cup_{1\leq i\leq r}X_{i}$, where $X_{i}$ is adjacent to exactly the vertices of $\otimes_{1\leq j\neq i\leq r}G^{j}$. Since $G$ is $F_{2}$-free, $G[X_{i}]$ has no vertex of degree at least 2 for any $1\leq i\leq r$. We must have $G\in\mathcal{G}_{n,r}$ as $G$ has maximum edges, as desired.

\medskip
 
 \f{\bf Case 2.} There is at most one $i$, such that $|H^{i}_{1}|\geq2$.
 
\medskip

 In this case, we can assume that $|H^{i}_{1}|=1$ for $i\geq2$ without loss of generality. That is $G_{i}=N\cdot K_{1}$ for $i\geq2$. Recall that $4m|N$.
 Also, let $V'=V(G)-V(\otimes_{1\leq i\leq r}G^{i})$. By symmetry, we can divide the vertices in $V'$ by $r+1$ parts as follows. Let $X_{1}$ be the set of vertices $u$ in $V'$, such that $u$ is adjacent to all the vertices of $\otimes_{2 \leq i\leq r}G^{i}$. For $2\leq i\leq r$, let $X_{i}$ be the set of vertices $u$ in $V'$, such that $u$ is adjacent to all the vertices of $\otimes_{1 \leq j\neq i\leq r}G^{j}$ but to no vertex of $V(G^{i})$. Let $Y=V'-\cup_{1\leq i\leq r}X_{i}$, i.e., the set of vertices which have non-neighbors in at least two among $G^{i}$ with $1\leq i\leq r$. Note that each vertex in $\cup_{2\leq i\leq r}X_{i}$ is adjacent to  $(r-1)N$ vertices of $\otimes_{1\leq i\leq r}G^{i}$. Each vertex in $Y$ is adjacent to at most $(r-1)N-1$ vertices of $\otimes_{1\leq i\leq r}G^{i}$ (using symmetry). The way connecting the vertices in $X_{1}$ to $V(G^{1})$ is not known now. 
 Note that
 \begin{equation}
\begin{aligned}
e(G)&=e(G[V'])+e(T(rN,r))+e(G^{1})+e_{G}(V(\otimes_{1\leq i\leq r}G^{i}),V')\\
&=e(G[V'])+e(T(rN,r))+e(G^{1})\\
&+e_{G}(V(\otimes_{1\leq i\leq r}G^{i}),X_{1})+e_{G}(V(\otimes_{1\leq i\leq r}G^{i}),\cup_{2\leq i\leq r}X_{i})+e_{G}(V(\otimes_{1\leq i\leq r}G^{i}),Y)\\
&\leq e(G[V'])+e(T(rN,r))+e(G^{1})\\
&+e_{G}(V(\otimes_{1\leq i\leq r}G^{i}),X_{1})+(r-1)N|\cup_{2\leq i\leq r}X_{i}|+((r-1)N-1)|Y|\\
&=e(G[V'])+e(T(rN,r))+e(G^{1})\\
&+e_{G}(V(\otimes_{1\leq i\leq r}G^{i}),X_{1})+(r-1)N\cdot(n-rN-|X_{1}|)-|Y|\\
&=e(G[V'])+e(T(rN,r))+e(G^{1})\\
&+e_{G}(V(G^{1}),X_{1})+(r-1)N|X_{1}|+(r-1)N\cdot(n-rN-|X_{1}|)-|Y|\\
&=e(G[V'])+e(G^{1})+e_{G}(V(G^{1}),X_{1})+e(T(rN,r))+(r-1)N\cdot(n-rN)-|Y|.
\end{aligned}\notag
\end{equation}
We will show that $G=K_{\frac{r}{2}}\otimes T(n-\frac{r}{2},r)$. Recall that $$e(G)\geq e(G')+e(T(rN,r))+\frac{r}{2}N+(n-rN)(r-1)N,$$
where $G'\in{\rm EX}(n-rN,\mathcal{F})$. Since the rest of the proof  is very similar to the Case 2 in $(i)$, we omit it.
This completes the proof.\hfill$\Box$

\medskip

\f{\bf Proof of Theorem \ref{m Icosahe}.} 
Let $$F^{1}=P_{6}\otimes 3K_{1}\otimes3K_{1},$$
$$F^{2}=(K_{1,2}\cup K_{2})\otimes(K_{2}\cup 2K_{1})\otimes3K_{1},$$
and $$F^{3}=2K_{3}\otimes3K_{1}\otimes3K_{1}.$$
Then, (see \cite{S2}) $I^{12}\subseteq F^{i}$ for any $1\leq i\leq3$. Thus, $mI^{12}\subseteq m F^{i}$ for any $1\leq i\leq3$. Recall that $\chi(I^{12})=4$ and $\chi(I^{12}-U)=4$ for any $U\subseteq V(I^{12})$ with $|U|=2$. Thus, $\chi(mI^{12}-U)=4$ for any $U\subseteq V(mI^{12})$ with $|U|=3m-1$. It follows that
$K_{3m-1}\otimes T(n-3m+1,3)$ is $mI^{12}$-free.
We will prove that ${\rm EX}(n,mI^{12})=\left\{K_{3m-1}\otimes T(n-3m+1,3)\right\}$ for sufficiently large $n$ using the progressive induction.

For $n\geq1$,  let $\mathfrak{U}_{n}={\rm EX}(n,mI^{12})$.
For a  graph $G'\in{\rm EX}(n,mI^{12})$, we say $G'$ has property $B$, if $G'=K_{3m-1}\otimes T(n-3m+1,3)$.  For a  graph $G\in{\rm EX}(n,mI^{12})$, define 
$$\phi(G)=e(G)-e(K_{3m-1}\otimes T(n-3m+1,3)).$$
 Clearly, $\phi(G)$ is a non-negative integer-valued function. If $G$ has property $B$, then $\phi(G)=0$. By Lemma \ref{progressive induction}, it suffices to show that for sufficiently large $n$, either $G$ satisfies property $B$, or there is a $\frac{n}{2}\leq n'<n$ and a $G'\in{\rm EX}(n',mI^{12})$ such that $\phi(G')>\phi(G)$.

Let $G\in {\rm EX}(n,mI^{12})$ with sufficiently large $n$. Recall that $mI^{12}\subseteq mF^{1},mF^{2}$. By Lemma \ref{Pl-free} (noting $t=12m$ there), there is a constant $N$ (by letting $N_{0}=12(6m-1)$ and  $N=O(r,t)\cdot N_{0}$ in Lemma \ref{Pl-free}), such that $G$  has an induced subgraph $\otimes_{1\leq i\leq 3}G^{i}$ satisfying:\\
 $(1)$  $|G^{i}|=N$ for any $1\leq i\leq 3$;\\
 $(2)$ for any $1\leq i\leq 3$, $G^{i}=\cup _{1\leq j\leq \ell_{i}}H^{i}_{j}$, where $H^{i}_{1},H^{i}_{2},...,H^{i}_{\ell_{i}}$ are symmetric subgraphs of $G$ and $\ell_{i}\geq 12(6m-1)$.\\
 Note that $12(6m-1)| N$.

 \medskip
 
 \f{\bf Case 1.} There are at least two $i$, such that $|H^{i}_{1}|\geq2$.
 
\medskip

In this case, we can assume that $|H^{i}_{1}|\geq2$ for $i=1,2$ without loss of generality. Note that $\ell_{i}\geq12m$. Since $G$ is $mF^{2}$-free, we must have that $|H^{1}_{1}|=|H^{2}_{1}|=2$ and $|H^{3}_{1}|\leq2$. That is, $G^{i}=\ell_{i}K_{1}$ or $G^{i}=\ell_{i}K_{2}$ for any $1\leq i\leq 3$. 
 
Let $V'=V(G)-V(\otimes_{1\leq i\leq 3}G^{i})$. Let $X$ be the set of vertices $u\in V'$ such that $u$ has at least $2N+1$ neighbors in $\otimes_{1\leq i\leq 3}G^{i}$. By symmetry, any vertex $u\in X$ has at most $\frac{N}{2}$ non-neighbors in $\otimes_{1\leq i\leq 3}G^{i}$. It is easy to see that $F^{2}$ is a subgraph of $G$ induced by $u$ and $\cup_{1\leq i\leq 3,1\leq j\leq6}V(H^{i}_{j})$. If $|X|\geq m$, assume that $u_{\ell}\in X$ with $1\leq \ell\leq m$. Then $mF^{2}$ is a subgraph of $G$ induced by $\left\{u_{1},u_{2},...,u_{m}\right\}$ and $\cup_{1\leq i\leq 3,1\leq j\leq6m}V(H^{i}_{j})$. This contradicts that $G$ is $mF^{2}$-free. Hence, $|X|\leq m-1$. Since $v\in V'-X$ is adjacent to at most $2N$ vertices of $\otimes_{1\leq i\leq 3}G^{i}$, we have 
 $$e_{G}(V(\otimes_{1\leq i\leq 3}G^{i}),V')\leq(n-3N)2N+|X|N\leq(n-3N)2N+(m-1)N.$$

By a calculation. we have 
$e(K_{3m-1}\otimes T(n-3m+1,3))=e(K_{3m-1}\otimes T(n-3N-3m+1,3))+e(T(3N,3))+(3m-1)N+(n-3N)2N.$
 Let $G'\in{\rm EX}(n-rN,\mathcal{F})$. Then $e(G[V'])\leq e(G')$. Note that $\phi(G')=e(G')-e(K_{3m-1}\otimes T(n-3N-3m+1,3))$.
 Thus,
\begin{equation}
\begin{aligned}
\phi(G)&=e(G)-e(K_{3m-1}\otimes T(n-3m+1,3))\\
&\leq e(G[V'])+e_{G}(V(\otimes_{1\leq i\leq 3}G^{i}),V')+e(T(3N,3))+\frac{3N}{2}\\
&-e(K_{3m-1}\otimes T(n-3N-3m+1,3))-e(T(3N,3))-(3m-1)N-(n-3N)2N\\
&\leq\phi(G')+(n-3N)2N+(m-1)N+\frac{3N}{2}-(3m-1)N-(n-3N)2N\\
&=\phi(G')-\frac{N}{2}(4m-3)\\
&<\phi(G'),
\end{aligned}\notag
\end{equation}
as desired.

\medskip
 
 \f{\bf Case 2.} There is at most one $i$, such that $|H^{i}_{1}|\geq2$.
 
\medskip

In this case, we can assume that $|H^{2}_{1}|=|H^{3}_{1}|=1$ without loss of generality. That is $G_{i}=N\cdot K_{1}$ for $i=2,3$. Recall that  $12(6m-1)| N$.
Let $V'=V(G)-V(\otimes_{1\leq i\leq 3}G^{i})$. By symmetry, we can divide the vertices in $V'$ by $4$ parts as follows. Let $X_{1}$ be the set of vertices $u$ in $V'$, such that $u$ is adjacent to all the vertices of $\otimes_{2 \leq i\leq 3}G^{i}$. For $2\leq i\leq 3$, let $X_{i}$ be the set of vertices $u$ in $V'$, such that $u$ is adjacent to all the vertices of $\otimes_{1 \leq j\neq i\leq 3}G^{j}$ but to no vertex of $V(G^{i})$. Let $Y=V'-\cup_{1\leq i\leq 3}X_{i}$. Note that each vertex in $\cup_{2\leq i\leq 3}X_{i}$ is adjacent to  $2N$ vertices of $\otimes_{1\leq i\leq 3}G^{i}$. Each vertex in $Y$ is adjacent to at most $2N-1$ vertices of $\otimes_{1\leq i\leq 3}G^{i}$ (using symmetry). 
 
Recall that $e(K_{3m-1}\otimes T(n-3m+1,3))=e(K_{3m-1}\otimes T(n-3N-3m+1,3))+e(T(3N,3))+(3m-1)N+(n-3N)2N$. Let $G'\in{\rm EX}(n-3N,mI^{12})$. Then $e(G[V'])\leq e(G')$. Note that $\phi(G')=e(G')-e(K_{3m-1}\otimes T(n-3N-3m+1,3))$. 
Thus,
 \begin{equation}
\begin{aligned}
\phi(G)&=e(G)-e(K_{3m-1}\otimes T(n-3m+1,2))\\
&=e(G[V'])+e_{G}(V(\otimes_{1\leq i\leq 3}G^{i}),V')+e(T(3N,3))+e(G^{1})\\
&-e(K_{3m-1}\otimes T(n-3N-3m+1,3))-e(T(3N,3))-(3m-1)N-(n-3N)2N\\
&\leq\phi(G')+e_{G}(V(\otimes_{1\leq i\leq 3}G^{i}),V')+e(G^{1})-(3m-1)N-(n-3N)2N.
\end{aligned}\notag
\end{equation}
If $e_{G}(V(\otimes_{1\leq i\leq 3}G^{i}),V')+e(G^{1})-(3m-1)N-(n-3N)2N<0$, then $\phi(G)<\phi(G')$, as desired. 
 Thus, we can assume that 
 $$e_{G}(V(\otimes_{1\leq i\leq 3}G^{i}),V')+e(G^{1})\geq (3m-1)N+(n-3N)2N.$$
  We will show that $G$ has property $B$ (i.e., $G=K_{3m-1}\otimes T(n-3m+1,3)$).

Note that
 \begin{equation}
\begin{aligned}
&e_{G}(V(\otimes_{1\leq i\leq 3}G^{i}),V')\\
&=e_{G}(V(\otimes_{1\leq i\leq 3}G^{i}),X_{1})+e_{G}(V(\otimes_{1\leq i\leq 3}G^{i}),\cup_{2\leq i\leq 3}X_{i})+e_{G}(V(\otimes_{1\leq i\leq 3}G^{i}),Y)\\
&\leq e_{G}(V(\otimes_{1\leq i\leq 3}G^{i}),X_{1})+2N|\cup_{2\leq i\leq 3}X_{i}|+(2N-1)|Y|\\
&=e_{G}(V(\otimes_{1\leq i\leq 3}G^{i}),X_{1})+2N\cdot(n-3N-|X_{1}|)-|Y|\\
&=e_{G}(V(G^{1}),X_{1})+2N|X_{1}|+2N\cdot(n-3N-|X_{1}|)-|Y|\\
&=e_{G}(V(G^{1}),X_{1})+2N\cdot(n-3N)-|Y|.
\end{aligned}\notag
\end{equation}
Recall $$e_{G}(V(\otimes_{1\leq i\leq 3}G^{i}),V')+e(G^{1})\geq (3m-1)N+(n-3N)2N.$$
It follows that 
$$e(G^{1})+e_{G}(V(G^{1}),X_{1})\geq(3m-1)N+|Y|\geq(3m-1)N.$$
Note that each vertex in $V(G^{1})\cup X_{1}$ is adjacent to all the vertices in $V(\otimes_{2\leq i\leq 3}G^{i})$. Thus, $G[V(G^{1})\cup X_{1}]$ is $P_{6m}$-free, since $G$ is $mF^{1}$-free. Note that $H^{1}_{1},H^{1}_{1},...,H^{1}_{\ell_{1}}$ are also symmetric in $G[V(G^{1})\cup X_{1}]$ and $|G^{1}|=N$. By Lemma \ref{desired symmetric subgraphs}, we have
$e(G^{1})+e_{G}(V(G^{1}),X_{1})=(3m-1)N$ (implying $|Y|=0$).
Moreover, either $Q^{1}_{i}=K_{6m-1}$ for any $1\leq i\leq\ell_{1}$ and $e_{G}(V(G^{1}),X_{1})=0$, or $Q^{1}_{i}=K_{1}$ for any $1\leq i\leq\ell_{1}$ and there are $3m-1$ vertices in $X_{1}$ which are all adjacent to the (unique) vertex in $Q^{1}_{i}$ for all $1\leq i\leq\ell_{1}$. If $Q^{1}_{i}=K_{6m-1}$ for any $1\leq i\leq\ell_{1}$, then $mF^{3}$ is contained in $\otimes_{1\leq i\leq 3}G^{i}$, a contradiction. Hence, we must have the latter case. That is, $G^{1}=N\cdot K_{1}$, and there are $3m-1$ vertices in $X_{1}$, say $y_{1},y_{2},...,y_{3m-1}$, such that they are all adjacent to all the vertices of $G^{1}$. Recall $Y=\emptyset$. Let $X=\left\{y_{1},y_{2},...,y_{3m-1}\right\}$. For convenience, we still use $X_{1}$ to denote $X_{1}-X$.
 For each $1\leq i\leq 3$, let $S_{i}=V(G^{i})\cup X_{i}$. Then $G$ has a partition $V(G)=X\cup(\cup_{1\leq i\leq 3}S_{i})$ satisfying:\\
$(3)$ $X=\left\{y_{1},y_{2},...,y_{3m-1}\right\}$ and $S_{i}=V(G^{i})\cup X_{i}$ for any $1\leq i\leq 3$, where $|V(G^{i})|=N\geq12m$;\\
$(4)$ for any $1\leq i\leq 3$, $N_{G}(v)=V(G)-S_{i}$ for any $v\in V(G^{i})$.

\medskip

For any $1\leq i\leq 3$, let $X'_{i}$ be the set of vertices in $S_{i}$ which have neighbors in $X$. Now we show that the vertices in $X'_{i}$ are isolated in $G[S_{i}]$. Suppose not. Without loss of generality, assume that $v\in X'_{1}$ such that $v$ is adjacent to $y_{3m-1}$ and adjacent to another vertex $w$ in $S_{1}$. Let $u_{1},u_{2},...,u_{3m-1}$ be other $3m-1$ vertices in $V(G^{1})$. Then $u_{1}y_{1}u_{2}y_{2}\cdots u_{3m-1}y_{3m-1}vw$ is a $P_{6m}$. Since the vertices of this path are all adjacent to the vertices in $V(\otimes_{2\leq i\leq 3}G^{i})$, we see $mF^{1}$ is contained in $G$, a contradiction. Hence, the vertices in $X'_{i}$ are isolated in $G[S_{i}]$ for any $1\leq i\leq 3$.

Note that $G[X_{i}-X'_{i}]$ is $P_{6m}$-free for any $1\leq i\leq 3$. By Lemma \ref{ex n Pt}, we have $e(G[X_{i}-X'_{i}])\leq(3m-1)|X_{i}-X'_{i}|$ with equality if and only if $G[X_{i}-X'_{i}]$ is the disjoint union of $K_{6m-1}$. 
Thus,
\begin{equation}
\begin{aligned}
e(G)&\leq e(K_{3m-1}\otimes K_{|S_{1}|,|S_{2}|,...,|S_{3}|})+(\sum_{1\leq i\leq 3}e(G[X_{i}-X'_{i}]))-(\sum_{1\leq i\leq 3}(3m-1)|X_{i}-X'_{i}|)\\
&\leq e(K_{3m-1}\otimes T(n-3m+1,3))+(\sum_{1\leq i\leq 3}e(G[X_{i}-X'_{i}]))-(\sum_{1\leq i\leq 3}(3m-1)|X_{i}-X'_{i}|)\\
&\leq e(K_{3m-1}\otimes T(n-3m+1,3)).
\end{aligned}\notag
\end{equation}
However, $e(G)\geq e(K_{3m-1}\otimes T(n-3m+1,3))$. This forces that $G[X]=K_{3m-1}$, and $||S_{i}|-\frac{n-3m+1}{3}|<1$ and $e(G[X_{i}-X'_{i}])=(3m-1)|X_{i}-X'_{i}|$ for any $1\leq i\leq 3$.
Then $G[X_{i}-X'_{i}]$ is the disjoint union of $K_{6m-1}$. If $X_{i}-X'_{i}$ is not empty for some $1\leq i\leq 3$, say $i=1$ without loss of generality, then $G[X_{1}-X'_{1}]$ contains a $K_{6m-1}$. Recall $G[X]=K_{3m-1}$. Then, $G[X\cup X_{1}]$ contains $2mK_{3}$ (as $m\geq2$).
Since all vertices of $\otimes_{2\leq i\leq 3}G^{i}$ are adjacent to vertices in $X\cup X_{1}$, $G$ contains a $mF^{3}$, a contradiction. Hence, $X_{i}-X'_{i}=\emptyset$ for any $1\leq i\leq 3$. This means that all the vertices in $S_{i}$ are isolated in $G[S_{i}]$ for any $1\leq i\leq 3$. Then, $G$ is a spanning subgraph of $K_{3m-1}\otimes T(n-3m+1,3)$. Hence, $G=K_{3m-1}\otimes T(n-3m+1,3)$, as desired.  This completes the proof.
\hfill$\Box$

\section{Proof of Theorem \ref{bi-matching-free}} 

In this section, we give the proof of Theorem \ref{bi-matching-free}.

\medskip

\f{\bf Proof of Theorem \ref{bi-matching-free}.} Recall that $\ell\geq3,m\geq2\ell$, and $\mathcal{F}=\left\{F_{1},F_{2},F_{3},F_{4}\right\}$, where
$$F_{1}=(P_{\ell}\cup\overline{K_{m}})\otimes T(m(r-1),r-1),$$
$$F_{2}=mK_{2}\otimes mK_{2}\otimes T(m(r-2),r-2),$$ 
$$F_{3}=(2K_{\ell-1}\cup\overline{K_{m}})\otimes (K_{2}\cup\overline{K_{m}})\otimes T(m(r-2),r-2),$$ 
$$F_{4}=(K_{\ell-1}\cup\overline{K_{m}})\otimes (K_{\ell-1}\cup\overline{K_{m}})\otimes T(m(r-2),r-2).$$ 
First, we can check that for any $1\leq i\leq r$, $F_{i}$ contains a $P_{\ell}$ after deleting any $(r-1)$ independents sets. Thus, $H\otimes T(n-|H|,r-1)$ is $\mathcal{F}$-free, if $H$ is $P_{\ell}$-free. For any $n\geq 2mr$, let $\mathcal{G}_{n}$ be the set of graphs of order $n$ and of form $H\otimes T(n-|H|,r-1)$, where $H$ is a $P_{\ell}$-free graph. Let $G_{n}$ an extremal graph with maximum edges in $\mathcal{G}_{n}$. By Lemma \ref{ex n Pt sharp}, we can always choose $G_{n}$ with $H=tK_{\ell-1}\cup K_{\ell_{0}}$, where $|H|=(\ell-1)t+\ell_{0}$ with $0\leq\ell_{0}<\ell-1$. Clearly, $0<|H|-\frac{n}{r}<r\ell$.

Let $G'$ be the  graph obtained from $G_{n}$ by deleting the vertices of $K_{\ell-1}$ in $H$, and deleting $\ell-1$ vertices in each part of $T(n-|H|,r-1)$ in $H$. Clearly, $G'\in\mathcal{G}_{n-r(\ell-1)}$ and $e(G')=e(G_{n})-(r-1)(\ell-1)(n-(\ell-1)r)-e(T((\ell-1)r,r))-\frac{\ell-2}{2}(\ell-1)$. Conversely, from $G_{n-r(\ell-1)}$ we can obtain a graph $G''\in\mathcal{G}_{n}$ by adding a $K_{\ell-1}$ in its first part, and adding $\ell-1$ independent vertices in each of the remained $r-1$ parts. Clearly, $e(G_{n-r(\ell-1)})=e(G'')-(r-1)(\ell-1)(n-(\ell-1)r)-e(T((\ell-1)r,r))-\frac{\ell-2}{2}(\ell-1)$.
Since $e(G')\leq e(G_{n-r(\ell-1)})$ and $e(G'')\leq e(G_{n})$, we have $e(G')=e(G_{n-r(\ell-1)})=e(G_{n})
-(r-1)(\ell-1)(n-(\ell-1)r)-e(T((\ell-1)r,r))-\frac{\ell-2}{2}(\ell-1)$. Similarly,  we can show
that for any constant $N$ with $(\ell-1)| N$, $e(G_{n})=e(G_{n-Nr})+e(T(Nr,r))+(n-rN)(r-1)N+\frac{\ell-2}{2}N$.

For $n\geq1$,  let $\mathfrak{U}_{n}={\rm EX}(n,\mathcal{F})$.
For a  graph $G'''\in{\rm EX}(n,\mathcal{F})$, we say $G'''$ has property $B$, if $G'''\in\mathcal{G}_{n}$.  For a  graph $G\in{\rm EX}(n,\mathcal{F})$, define 
$$\phi(G)=e(G)-e(G_{n}).$$
 Clearly, $\phi(G)$ is a non-negative integer-valued function. If $G$ has property $B$, then $\phi(G)=0$. By Lemma \ref{progressive induction}, it suffices to show that for sufficiently large $n$, either $G$ satisfies property $B$, or there is a $\frac{n}{2}\leq n'<n$ and a $G_{0}\in{\rm EX}(n',\mathcal{F})$ such that $\phi(G_{0})>\phi(G)$.

Since $F_{1},F_{2}\in \mathcal{F}$, we see that $\mathcal{F}$ satisfies the condition in Lemma \ref{Pl-free}. Since $G\in {\rm EX}(n,\mathcal{F})$ with sufficiently large $n$, there is a constant $N$ (by letting $N_{0}=3m(\ell-1)$ and  $N=O(r,t)\cdot N_{0}$ in Lemma \ref{Pl-free}), such that $G$  has an induced subgraph $\otimes_{1\leq i\leq r}G^{i}$ satisfying:\\
 $(1)$  $|G^{i}|=N$ for any $1\leq i\leq r$;\\
 $(2)$ for any $1\leq i\leq r$, $G^{i}=\cup _{1\leq j\leq \ell_{i}}H^{i}_{j}$, where $H^{i}_{1},H^{i}_{2},...,H^{i}_{\ell_{i}}$ are symmetric subgraphs of $G$ and $\ell_{i}\geq 3m$.\\
 Note that $3m(\ell-1)| N$.  
  Since $G$ is $F_{2}$-free, there is at most one $i$, such that $|H^{i}_{1}|\geq2$. Without loss of generality, we can assume that $|H^{i}_{1}|=1$ for $i\geq2$. That is $G_{i}=N\cdot K_{1}$ for $i\geq2$. Let $V'=V(G)-V(\otimes_{1\leq i\leq r}G^{i})$. By symmetry, we can divide the vertices in $V'$ by $r+1$ parts as follows. Let $X_{1}$ be the set of vertices $u$ in $V'$, such that $u$ is adjacent to all the vertices of $\otimes_{2 \leq i\leq r}G^{i}$. For $2\leq i\leq r$, let $X_{i}$ be the set of vertices $u$ in $V'$, such that $u$ is adjacent to all the vertices of $\otimes_{1 \leq j\neq i\leq r}G^{j}$ but to no vertex of $V(G^{i})$. Let $Y=V'-\cup_{1\leq i\leq r}X_{i}$, i.e., the set of vertices which have non-neighbors in at least two among $G^{i}$ with $1\leq i\leq r$. Note that each vertex in $\cup_{2\leq i\leq r}X_{i}$ is adjacent to  $(r-1)N$ vertices of $\otimes_{1\leq i\leq r}G^{i}$. Each vertex in $Y$ is adjacent to at most $(r-1)N-1$ vertices of $\otimes_{1\leq i\leq r}G^{i}$ (using symmetry). The way connecting the vertices in $X_{1}$ to $V(G^{1})$ is not known now. 
 
Recall that $e(G_{n})=e(G_{n-Nr})+e(T(Nr,r))+(n-rN)(r-1)N+\frac{\ell-2}{2}N$.  Let $G_{0}\in{\rm EX}(n-rN,\mathcal{F})$. Then $e(G[V'])\leq e(G_{0})$. Note that $\phi(G_{0})=e(G_{0})-e(G_{n-rN})$. 
Thus,
 \begin{equation}
\begin{aligned}
\phi(G)&=e(G)-e(G_{n})\\
&=e(G[V'])+e_{G}(V(\otimes_{1\leq i\leq r}G^{i}),V')+e(T(rN,r))+e(G^{1})\\
&-(e(G_{n-Nr})+e(T(Nr,r))+(n-rN)(r-1)N+\frac{\ell-2}{2}N)\\
&\leq\phi(G_{0})+e_{G}(V(\otimes_{1\leq i\leq r}G^{i}),V')+e(G^{1})-\frac{\ell-2}{2}N-(n-rN)(r-1)N.
\end{aligned}\notag
\end{equation}
If $e_{G}(V(\otimes_{1\leq i\leq r}G^{i}),V')+e(G^{1})-\frac{\ell-2}{2}N-(n-rN)(r-1)N<0$, then $\phi(G)<\phi(G_{0})$, as desired. 
 Thus, we can assume that 
 $$e_{G}(V(\otimes_{1\leq i\leq r}G^{i}),V')+e(G^{1})\geq\frac{\ell-2}{2}N+(n-rN)(r-1)N.$$
  We will show that $G$ has property $B$ (i.e., $G\in\mathcal{G}_{n}$).

Similar to Theorem \ref{general Icosahe}, we can show that
 $e_{G}(V(\otimes_{1\leq i\leq r}G^{i}),V')=e_{G}(V(G^{1}),X_{1})+(r-1)N\cdot(n-rN)-|Y|$.
Recall $e_{G}(V(\otimes_{1\leq i\leq r}G^{i}),V')+e(G^{1})\geq\frac{\ell-2}{2}N+(n-rN)(r-1)N$.
It follows that 
$$e(G^{1})+e_{G}(V(G^{1}),X_{1})\geq\frac{\ell-2}{2}N+|Y|\geq\frac{\ell-2}{2}N.$$
Since $G$ is $F_{1}$-free, we have $G[X_{1}\cup V(G^{1})]$ is $P_{\ell}$-free.
By Lemma \ref{desired symmetric subgraphs}, we have
$e(G^{1})+e_{G}(V(G^{1}),X_{1})=\frac{\ell-2}{2}N$ (implying $|Y|=0$).
Moreover, either $Q^{1}_{i}=K_{\ell-1}$ for any $1\leq i\leq\ell_{1}$ and $e_{G}(V(G^{1}),X_{1})=0$, or ($\ell$ is even) $Q^{1}_{i}=K_{1}$ for any $1\leq i\leq\ell_{1}$ and there are $\frac{\ell-2}{2}$ vertices in $X_{1}$ which are all adjacent to the (unique) vertex in $Q^{1}_{i}$ for all $1\leq i\leq\ell_{1}$. Recall that $Y=\emptyset$. For any $1\leq i\leq r$, each vertex in $V(Q^{i})$ is adjacent to all the vertices in $\cup_{1\leq j\neq i\leq r}(X_{j}\cup V(G^{j}))$. If $Q^{1}_{i}=K_{\ell-1}$ for any $1\leq i\leq\ell_{1}$, then $G[X_{i}\cup V(G^{i})]$ has no edges for any $2\leq i\leq r$ (otherwise $F_{3}$ will be contained in $G$, a contradiction). Thus, $G\in\mathcal{G}_{n}$, as desired. It remains to consider the latter case. That is, ($\ell$ is even) $G^{1}=N\cdot K_{1}$, and there are $\frac{\ell-2}{2}$ vertices in $X_{1}$, say $y_{1},y_{2},...,y_{\frac{\ell-2}{2}}$, such that they are all adjacent to all the vertices of $G^{1}$. Let $X=\left\{y_{1},y_{2},...,y_{\frac{\ell-2}{2}}\right\}$. For convenience, we still use $X_{1}$ to denote $X_{1}-X$.
 For each $1\leq i\leq r$, let $S_{i}=V(G^{i})\cup X_{i}$. Then $G$ has a partition $V(G)=X\cup(\cup_{1\leq i\leq r}S_{i})$ satisfying:\\
$(3)$ $X=\left\{y_{1},y_{2},...,y_{\frac{\ell-2}{2}}\right\}$ and $S_{i}=V(G^{i})\cup X_{i}$ for any $1\leq i\leq r$, where $|V(G^{i})|=N\geq3m$;\\
$(4)$ for any $1\leq i\leq r$, $N_{G}(v)=V(G)-S_{i}$ for any $v\in V(G^{i})$.

For any $1\leq i\leq r$, let $X'_{i}$ be the set of vertices in $X_{i}$ which have neighbors in $X$. Similar to Theorem \ref{general Icosahe}, we can show that the vertices in $X'_{i}$ are isolated in $G[S_{i}]$ for any $1\leq i\leq r$. Recall that $G[X_{i}-X'_{i}]$ is $P_{\ell}$-free for any $1\leq i\leq r$. By Lemma \ref{ex n Pt}, we have $e(G[X_{i}-X'_{i}])\leq\frac{\ell-2}{2}|X_{i}-X'_{i}|$ with equality if and only if $G[X_{i}-X'_{i}]$ is the disjoint union of several $K_{\ell-1}$. Let $\overline{e}$ be the number of non-edges between two distinct $S_{i}$ and $S_{j}$.
Thus,
\begin{equation}
\begin{aligned}
e(G)&\leq e(K_{\frac{\ell-2}{2}}\otimes K_{|S_{1}|,|S_{2}|,...,|S_{r}|})+(\sum_{1\leq i\leq r}e(G[X_{i}-X'_{i}]))-(\sum_{1\leq i\leq r}\frac{\ell-2}{2}|X_{i}-X'_{i}|)-\overline{e}\\
&\leq e(K_{\frac{\ell-2}{2}}\otimes T(n-\frac{\ell-2}{2},r))+(\sum_{1\leq i\leq r}e(G[X_{i}-X'_{i}]))-(\sum_{1\leq i\leq r}\frac{\ell-2}{2}|X_{i}-X'_{i}|)-\overline{e}\\
&\leq e(K_{\frac{\ell-2}{2}}\otimes T(n-\frac{\ell-2}{2},r))-\overline{e}.
\end{aligned}\notag
\end{equation}
However, $e(G)\geq e(K_{\frac{\ell-2}{2}}\otimes T(n-\frac{\ell-2}{2},r))$ as $G\in{\rm EX}(n,\mathcal{F})$. This forces that $\overline{e}=0$  and $e(G[X_{i}-X'_{i}])=\frac{\ell-2}{2}|X_{i}-X'_{i}|$ for any $1\leq i\leq r$.
Then $G[X_{i}-X'_{i}]$ is the disjoint union of several $K_{\ell-1}$. If $X_{i}-X'_{i}$ is not empty for two $i$, say $i=1,2$ without loss of generality, then both $G[X_{1}-X'_{1}]$ and $G[X_{2}-X'_{2}]$ contains a $K_{\ell-1}$. Since $\overline{e}=0$, there are no non-edges between these two $K_{\ell-1}$. Thus, $F_{4}$ is contained in $G$, a contradiction. Hence, $X_{i}-X'_{i}\neq\emptyset$ for at most one $i$. Without loss of generality, we can assume that  $X_{i}-X'_{i}=\emptyset$ for any $2\leq i\leq r$. This means that all the vertices in $S_{i}$ are isolated in $G[S_{i}]$ for any $2\leq i\leq r$. Note that $G[X\cup S_{1}]$ is $P_{\ell}$-free, since $G$ is $F_{1}$-free. Then, $G\in\mathcal{G}_{n}$, as desired.
This completes the proof.\hfill$\Box$

\medskip

\f{\bf Declaration of competing interest}

\medskip

There is no conflict of interest.

\medskip

\f{\bf Data availability statement}

\medskip

No data was used for the research described in the article.

\end{document}